\documentclass[a4paper]{scrartcl}
\usepackage[a4paper, total={6.5in, 9in}]{geometry}
\usepackage[utf8]{inputenc}
\usepackage{amsmath,amssymb,amsthm}
\usepackage{graphicx}
\usepackage{bm}
\usepackage{mathrsfs}
\usepackage{xcolor}
\usepackage{tensor}
\usepackage{hyperref}
\usepackage{caption}
\usepackage{colonequals}
\usepackage{cleveref}
\usepackage{authblk}
\usepackage{mathtools}
\hypersetup{colorlinks,
  linkcolor=blue,
  citecolor=blue,
  urlcolor=blue
}

\usepackage{comment}
\newcommand\blfootnote[1]{%
  \begingroup
  \renewcommand\thefootnote{}\footnote{#1}%
  \addtocounter{footnote}{-1}%
  \endgroup
}

\newtheorem{problem}{Problem}

\newtheorem{remark}{Remark}

\newcommand{\M}{\mathcal{M}}                
\newcommand{\R}{\mathbb{R}}                 

\newcommand{\Domain}{\Omega}
\newcommand{\RefDomain}{\hat{\Domain}}
\newcommand{\Surface}{\M}                   
\newcommand{\TangentBundle}{\mathrm{T}}     
\newcommand{\TensorBundle}{\mathrm{T}}      

\newcommand{\Basis}{\bm{t}}                 
\newcommand{\BasisOrth}{\tilde{\Basis}}     
\newcommand{\Metric}{\mathbf{g}}            
\newcommand{\MetricOrth}{\tilde{\Metric}}   
\newcommand{\MetricOrthCoeff}[1]{h_{(#1)}}  
\newcommand{\ChristSymb}[2]{\Gamma_{#1}^{#2}}  
\newcommand{\ChartMap}{\bm{\mu}}            
\newcommand{\GaussCurvature}{K}             

\newcommand{\Tri}[1]{\mathscr{S}_{#1}}      
\newcommand{\Elem}{S}                       
\newcommand{\dt}{\Delta{t}}                 

\newcommand{\Avg}[1]{\langle{#1}\rangle}    

\newcommand{\Emb}[1]
    {\mkern 1.5mu\underline{\mkern-1.5mu#1\mkern-1.5mu}\mkern 1.5mu}
\newcommand{\Ext}[1]
    {\mkern 1.5mu\overline{\mkern-1.5mu#1\mkern-1.5mu}\mkern 1.5mu}
\newcommand{\TT}[1]{\Emb{\bm{#1}}}          

\newcommand{\Vh}{V}                         
\newcommand{\VVh}{\bm{V}}                   

\renewcommand{\P}{\bm{P}}                   
\newcommand{\PP}{\mathcal{P}}               
\newcommand{\Ph}{\P_h}                      
\newcommand{\PPh}{\PP_h}                    
\newcommand{\QQ}{\mathcal{Q}}               
\newcommand{\QQh}{\QQ_h}                    

\newcommand{\grad}{\nabla}
\newcommand{\Grad}{\bm{\nabla}}

\newcommand{\Laplace}{\ensuremath{\bm{\Delta}}}

\newcommand{\Div}{\operatorname{div}}
\newcommand{\D}[2]{\frac{\partial{#1}}{\partial{#2}}}

\newcommand{\scalarProd}[2]{\ensuremath{\big\langle{#1}\,,\,{#2}\big\rangle}}

\newcommand{\Inner}[2]{\left({#1}\,,\,{#2}\right)}

\newcommand{\Abs}[1]{\left\lvert{#1}\right\rvert}
\newcommand{\Norm}[1]{\lVert{#1}\rVert}

\newcommand{\Set}[1]{\big\lbrace{#1}\big\rbrace}
\newcommand{\EnergyNorm}[1]{{\lvert\kern-0.25ex\lvert\kern-0.25ex\lvert{#1}\rvert\kern-0.25ex\rvert\kern-0.25ex\rvert}}

\newcommand{\Tdot}[1]{%
    \underset{{\raisebox{.4ex}[0pt][0pt]{\ensuremath{\scriptstyle{#1}}}}}{\cdot}}

\newcommand{\mb}{\bm{m}}
\newcommand{\nb}{\bm{n}}
\newcommand{\pb}{\bm{p}}

\newcommand{\ub}{\bm{u}}
\newcommand{\vb}{\bm{v}}

\newcommand{\xb}{\bm{x}}
\newcommand{\yb}{\bm{y}}

\newcommand{\sv}{\hat{\bm{s}}}
\newcommand{\xv}{\hat{\bm{x}}}

\usepackage{tikz}
\newcommand{\doi}[1]{doi:\href{https://doi.org/#1}{\detokenize{#1}}}

\date{\vspace*{-2cm}}

\title{Diffusion of tangential tensor fields: numerical issues and influence of geometric properties}

\author[1]{E. Bachini}
\author[2]{P. Brandner}
\author[2]{T. Jankuhn}
\author[1]{M. Nestler}
\author[1]{\mbox{S. Praetorius}\footnote{Corresponding author: \url{simon.praetorius@tu-dresden.de}}}
\author[2]{\mbox{A. Reusken}}
\author[1]{A. Voigt}
\affil[1]{Institut für Wissenschaftliches Rechnen, Technische Universität Dresden, D-01062 Dresden, Germany.}
\affil[2]{Institut für Geometrie und Praktische Mathematik, RWTH Aachen University, D-52056 Aachen, Germany}
\begin{document}
\maketitle
\blfootnote{\today}

\begin{abstract}
\noindent\textbf{Abstract}.\space%
We study the diffusion of tangential tensor-valued data on curved surfaces. For this purpose, several finite-element-based numerical methods are collected and used to solve a tangential surface $n$-tensor heat flow problem. These methods differ with respect to the surface representation used, the geometric information required, and the treatment of the tangentiality condition. We emphasize the importance of geometric properties and their increasing influence as the tensorial degree changes from $n=0$ to $n \geq 1$. A specific example is presented that illustrates how curvature drastically affects the behavior of the solution.

\vspace{1em}
\end{abstract}

\section{Introduction}
The Isaac Newton Institute program ``Computational Challenges in Partial Differential Equations'' in 2003 stimulated intensive research on surface partial differential equations (PDEs) in the mathematical community covering topics in modeling, numerical analysis, and applications. PDEs defined on curved surfaces are inherently nonlinear and require a geometric framework. An important breakthrough in the development of numerical methods for this type of PDEs is the avoidance of charts and atlases. Most commonly used methods are either based on a triangulated surface and require geometric information  through knowledge of the vertices and discrete normals, or are based on a level set technique where the geometric information is derived from the level set function. Most of these works deal with scalar-valued surface PDEs, see \cite{DE2013Finite,OlshanskiiReusken2017,BD2020Finite} for reviews of such finite-element-based approaches.
In the scalar case, the coupling between the surface geometry and the PDE solution is relatively weak, and  numerical approaches developed for PDEs in a flat space need only minor modifications to be applicable to surface equations, see, e.g., \cite{DE2013Finite,VV2006AMDiS}.
For $n$-tensor-valued surface PDEs with $n \geq 1$, these approaches are not directly applicable. The  tensor-fields must be considered as elements of the tangent bundle and the surface derivatives require more geometric information. This leads to a stronger influence of the surface geometry on the solution of the PDE.

In this paper we investigate this change in numerical complexity when moving from PDEs in flat domains to curved surfaces. We consider a specific class of problems, namely that of a surface heat equation for $n$-tensor fields on a smooth curved surface embedded in $\R^3$. We will focus on tensor ranks $n=0,1,2$. In the rest of the paper we will use for $n$-tensors, $n=0,1,2$, the terminology scalars, vectors and tensors, respectively. For $n=1,2$ the solution must be tangential. The numerics are restricted to finite element discretizations in space combined with low order BDF time stepping schemes, cf. \cite{LMV2013Backward}. When moving from PDEs in flat Euclidean domains to PDEs on a curved surface the following additional numerical issues arise:
\begin{itemize}
  \item[a)] \emph{Surface representation}. In flat domains, only the boundary of the computational domain needs to be represented or approximated. For problems on curved surfaces, however, the entire domain must be prescribed or discretized. An issue directly related to this is the quadrature used in the finite element method.
  \item[b)] \emph{Representation of the gradient operator and  geometry information}. For the  (covariant) surface gradient operator, different natural representations are available, which lead to different numerical approaches. The discretization process requires approximations of geometric quantities such as surface normals and curvature. We will see that the required geometric information depends on the representation of the gradient operator used and on the tensorial degree $n$.
  \item[c)] \emph{Tangentiality condition}. For $n$-tensor fields with $n \geq 1$ one has to take into account the condition that the solution must be tangential.
\end{itemize}
In recent years, several approaches have been developed to deal with these problems, leading to different numerical discretization methods. \emph{We present, in a unified framework, four methods known from the literature and explicitly address the different approaches these methods take with respect to a)--c).}
These four methods are:  A surface finite element method (SFEM) \cite{NNPV2018Orientational,NN2019Finite,HL2020Analysis,HaPr2021Tangential}, which extends the SFEM for scalar-valued surface PDEs \cite{DE2013Finite,Demlow2009HigherOrder} to tensor-valued surface PDEs; an intrinsic surface finite element method (ISFEM), which so far has only been considered for scalar-valued surface PDEs \cite{BFP2021Intrinsic}; a trace finite element method (TraceFEM) \cite{JankuhnReusken2020}, which extends the scalar version \cite{OlshanskiiReusken2017} to vector-valued PDEs; and a diffuse interface approach (DI) \cite{NNPV2018Orientational,NV2023diffuse}, which extends the approach for scalar-valued PDEs \cite{RV2006PDEs}.
We note that only very few rigorous discretization error analyses are available for vector- or tensor-valued surface PDEs. Such analyses for SFEM and TraceFEM applied to a vector-Laplace problem are given in \cite{HL2020Analysis,JankuhnReusken2020,HaPr2021Tangential}.

One conclusion from this comparative study is that for all four methods there is \emph{a significant increase in numerical complexity when moving from the scalar case to the vector- or tensor-valued problem}, which goes far beyond the increase in complexity in flat Euclidean domains. Depending on the geometry, an approximation of geometric properties that is sufficient to achieve the desired accuracy of the solution for the scalar case may fail for the vector or tensor case, cf. \Cref{sec:Comparison} for a further discussion.

We also consider the \emph{influence of the geometry on the solution of an $n$-tensor heat flow problem}. This is done on a surface with a rather simple geometry. We present results of numerical simulations using the four methods, which show that curvature drastically affects the behavior of the solution.

We would like to mention other vector- and matrix-valued finite element methods that are available in the literature: The $H(\operatorname{div})$ and $H(\operatorname{curl})$ conforming Brezzi--Douglas--Marini/Raviart--Thomas and N\'ed\'elec elements for vector fields, and Hellan--Herrmann--Johnson and Regge elements for tensor fields, see \cite{LLS2020Divergence,BDL2020Divergence,NeunteufelSchoeberl2019Hellan,Gawlik2020High}. These methods directly consider the tangent space of the surface by using the Piola or the convariant transformations, but require more sophisticated discretization approaches such as hybrid discontinuous Galerkin or mixed formulations. For these reasons, they are not considered in this paper.

The paper is structured as follows: In \Cref{sec:2}, we recall different surface representations. We also discuss different ways of representing tensors and gradient operators. Furthermore, we introduce the surface $n$-tensor-valued heat equation and summarize known analytical results. In \Cref{sec:FiniteElementDiscretizationSchemes}, we briefly describe the four numerical methods and discuss the above mentioned numerical issues a)--c). In \Cref{sec:NumericalExperiments}, the $n$-tensor-valued heat equation on a specific surface is solved numerically. Certain influences of the geometry on the behavior of the solution are discussed in the Sections~\ref{sec:Phenomenon1} -- \ref{sec:Phenomenon3}.

In the following we restrict ourselves to tensorial degree $n \leq 2$. This restriction is not essential for the results presented or for the applicability of the numerical methods, but it allows a clearer presentation. We provide reference solutions that can serve as benchmark problems. 

\section{Surface tensor diffusion}\label{sec:2}
\subsection{Surface representation}
Let $\Surface$ be a compact, orientable, two-dimensional surface isometrically embedded in $\R^3$. We consider two representations of this surface, namely based on a local  parametrization and as the zero level of a level set function. The tangent bundle of $\Surface$ is denoted by $\TangentBundle\Surface$ and for each $\xb\in\Surface$ a normal vector $\nb(\xb)\in\R^3$ is defined as the unit vector orthogonal to all tangent vectors in $\TangentBundle_{\xb}\Surface$.

\subsubsection{Parameterized surface}
\label{subsec:parametrizedsurface}
We assume that $\Surface$ can be covered by a $C^k$-atlas $\Set{(\ChartMap_r,\RefDomain_r,U_r)}_r$ of bijective mappings $\ChartMap_r\colon\RefDomain_r\to U_r\cap\Surface$, which are parametrizations of class $C^k$ with the domain open subsets $\RefDomain_r\subset\R^2$. We further assume that the transition maps $\ChartMap_r^{-1}\circ\ChartMap_s$ for overlapping co-domains, $\ChartMap_r(\RefDomain_r)\cap\ChartMap_s(\RefDomain_s)\neq\emptyset$, are $C^k$-diffeomorphisms.

For a local parametrization $\ChartMap = \ChartMap_r$, we denote the surface coordinate by $\xb=\ChartMap(\xv)\in\Surface$, with $\xv=(\hat{x}^1,\hat{x}^2)\in\RefDomain=\RefDomain_r$, the Jacobian of the parametrization by
\[
  [\bm{J}(\xb)]\indices{^\alpha_j} \colonequals \frac{\partial \mu^\alpha}{\partial\hat{x}^j}(\xv),\quad \alpha=1,2,3; j=1,2,
\]
and the surface metric tensor by $\Metric=\bm{J}^T\bm{J}$, i.e., $[\Metric]_{ij}=J\indices{^\alpha_i}J\indices{^\beta_j}\delta_{\alpha\beta}$. Here and in the rest of the text we use the Ricci calculus and the summation convention with Greek indices for the Cartesian coordinates in $\R^3$ and Latin indices for the coordinates in the parameter domain $\hat{\Omega}\subset\R^2$. Bold symbols refer to the vector or tensor object, while light symbols with indices refer to its components.

The columns of $\bm{J}$ are $\R^3$-vectors tangent to $\Surface$, i.e., $[\Basis_j(\xb)]^\alpha\colonequals J(\xb)\indices{^\alpha_j}$, $\Basis_j(\xb)\in\TangentBundle_{\xb}\Surface$ for $\alpha=1,2,3; j=1,2$, with $\xb=\ChartMap(\xv)$. This results in the definition of the normal direction field $\mb=\Basis_1\times\Basis_2$ and the corresponding unit normal field $\nb=\mb/\Norm{\mb}$.

A regular $C^2$-surface has an invertible metric. This allows to transform derivatives from the parameter domain $\RefDomain$ into surface derivatives, cf. \Cref{secintri}. The Weingarten map is given by
\[
  [\bm{H}(\xb)]\indices{^{\alpha\beta}} \colonequals -J\indices{^\beta_i}\,g^{ij}\frac{\partial n^\alpha}{\partial\hat{x}^j}(\xv)\,,
\]
where $\Metric^{-1}=[g^{ij}]$ is the inverse of the metric tensor.

\subsubsection{Level set characterization of the surface}
An implicit representation of $\Surface$ can be based on a $C^k$-mapping $\phi\colon\Domain\to\R$, where $\Surface\subset\Domain\subset\R^3$ is a three-dimensional domain containing the surface. We assume that $\grad\phi\neq 0$ on $\Surface$ and represent the surface as the zero-level set of $\phi$:
\[
  \Surface = \Set{\xb\in\Domain\mid\phi(\xb)=0}\,.
\]
In a sufficiently small $\delta$-neighborhood $U_\delta(\Surface)\subset\R^3$ of $\Surface$, we can define the normal direction $\Ext{\mb}(\xb)=\grad\phi(\xb)$, $\xb \in U_\delta(\Surface)$, and the normal field $\Ext{\nb}_\phi=\Ext{\mb}/\Norm{\Ext{\mb}}$ with $\Ext{\nb}_\phi\vert_{\Surface}=\nb$. Here and in the rest of the text we use an overline notation, e.g., $\Ext{\mb}$, to denote quantities that are defined not only on $\Surface$, but in a (small) three-dimensional neighborhood of $\Surface$.

A natural choice for $\phi$ would be the signed-distance function $\rho(\xb)=\operatorname{dist}(\xb,\Surface)$ with the property $\Norm{\grad\rho}\equiv 1$. Let $\delta>0$ be sufficiently small so that the closest-point projection $\pi\colon U_\delta(\Surface)\to\Surface$  is uniquely defined implicitly by
\begin{equation}\label{eq:closest-point-projection}
  \pi(\xb) = \xb - \rho(\xb)\nb(\pi(\xb)),\quad\xb\in U_\delta(\Surface)\,.
\end{equation}
Using the closest-point projection, the signed-distance function can be determined based on $\rho(\xb) = (\xb - \pi(\xb))\cdot\nb(\pi(\xb))$. The normal field $\Ext{\nb}(\xb)=\grad\rho(\xb)$, $\xb\in U_\delta(\Surface)$,   is a constant extension of the surface normal, i.e., $\Ext{\nb}(\xb)=\nb(\pi(\xb))$.
An alternative representation of the extended Weingarten map is given by $\Ext{\bm{H}}(\xb)=-\grad\Ext{\nb}(\xb)=-\grad^2\rho(\xb)$, for $\xb \in U_\delta(\Surface)$ with $\Ext{\bm{H}}\vert_{\Surface}=\bm{H}$.

\subsection{Representation of tensor fields and gradient operators}
\subsubsection{Intrinsic representation}%
\label{secintri}
Starting from the definition of a parameterized surface, one can  represent tensor fields and define derivatives using local coordinates in a reference domain $\hat{\Omega}$.

One way to choose the local coordinates is to consider the tangent vectors $\Basis_1,\Basis_2$, which are naturally associated with the parametrization $\ChartMap$, as the reference frame for the tangent plane $\TangentBundle_{\xb}\Surface$. We can then describe a function on $\Surface$ in the local coordinates and define the (intrinsic) surface gradients. Let $\ub^{(0)}\colon\Surface\to\R$ be a scalar differentiable function on $\Surface$, $\ub^{(1)}\colon\Surface \to \TangentBundle\Surface$ be a tangent vector field given by $\ub^{(1)}=u^1\Basis_{1} + u^2\Basis_{2}=u^i\Basis_{i}$, and $\ub^{(2)}\colon\Surface\to  T^2 \Surface $  a tangent tensor field given by $\ub^{(2)}=u^{ij} \Basis_{i} \otimes \Basis_{j}$. At a point $\xb\in\Surface$, the tensors of the contravariant components are denoted by underline notation, i.e, $\Emb{\ub}^{(0)}=u\in\R$, $\Emb{\ub}^{(1)}=[u^i] \in \R^2$, and $\Emb{\ub}^{(2)}=[u^{ij}]\in \R^{2 \times 2}$.
The intrinsic gradient of an $n$-tensor field is an $(n+1)$-tensor field and is defined by the following expressions in terms of the contravariant components:
\begin{align}
  \left[\Grad_{\Surface}\ub^{(0)}\right]^i &= g^{il} \D{u}{\hat{x}^l}, \label{defG1} \\
  \left[\Grad_{\Surface}\ub^{(1)}\right]^{ij} &=
  g^{il}\Grad_l u^{j} =
  g^{il}\left(\D{u^{j}}{\hat{x}^l}
  + \ChristSymb{lk}{j}u^{k}\right), \label{defG2} \\
  \left[\Grad_{\Surface}\ub^{(2)}\right]^{ijk} &=
  g^{il}\Grad_l u^{jk} =
  g^{il}\left(\D{u^{jk}}{\hat{x}^l}
  + \ChristSymb{lh}{j}u^{hk}
  + \ChristSymb{lh}{k}u^{jh}\right), \label{defG3}
\end{align}
where $\ChristSymb{ij}{k}=g^{kl}\Basis_l \cdot \tfrac{\partial \Basis_i}{\partial \hat x^j}$ denote the Christoffel symbols, for $i,j,k=1,2$. Note that these contain curvature information.
Scalar products of tensors are explicitly written in terms of the metric $\Metric$, e.g.,
\[
  \scalarProd{\TT{u}^{(1)}}{\TT{v}^{(1)}}_{\Metric} \colonequals u^i\,v_i = g_{ij}\,u^i\,v^{j}\,,
  \quad
  \scalarProd{\TT{u}^{(2)}}{\TT{v}^{(2)}}_{\Metric} \colonequals u^{ij}\,v_{ij} = g_{il}\,g_{jm}\,u^{ij}\,v^{lm}\,.
\]

To simplify the computation and to increase numerical stability, we also consider an orthogonal reference frame as the basis for the tangent plane $\TangentBundle_{\xb}\Surface$. For this we orthogonalize the vector $\Basis_{2}$ with respect to $\Basis_{1}$. This orthogonalization results in the orthogonal frame $\BasisOrth_{1},\BasisOrth_{2}$ on $\TangentBundle_{\xb}\Surface$, associated with the local coordinates $\sv=(\hat{s}^1,\hat{s}^2)$.
The corresponding  metric tensor is given by
\[
  \MetricOrth \colonequals \left(
    \begin{array}{ccc}
      \Norm{\BasisOrth_{1}}^2 & 0 \\
      0 &\Norm{\BasisOrth_{2}}^2  \\
    \end{array}
  \right) \equalscolon \left(
    \begin{array}{ccc}
      \MetricOrthCoeff{1}^2&     0         \\
      0             &\MetricOrthCoeff{2}^2 \\
    \end{array}
  \right)\,.
\]
We can now write explicit expressions for the Christoffel symbols in the basis $\{\BasisOrth_{1},\BasisOrth_{2}\}$:
\begin{align}\label{eq:crhistsymb-orth}
  \ChristSymb{ik}{k}=\ChristSymb{ki}{k}=\frac{1}{\MetricOrthCoeff{k}}\D{\MetricOrthCoeff{k}}{\hat{s}^i}\quad& i,k=1,2\,,
  \qquad\quad
  \ChristSymb{ii}{k}=-\frac{\MetricOrthCoeff{i}}{\MetricOrthCoeff{k}^2}\D{\MetricOrthCoeff{i}}{\hat{s}^k}\quad i\ne k\,,\\
  &\ChristSymb{ij}{k}=0\quad i\ne j\ne k\,,
\end{align}
which are simplified due to the orthogonality property.
Also the scalar products simplify with the metric $\MetricOrth$, e.g., $\scalarProd{\TT{u}^{(2)}}{\TT{v}^{(2)}}_{\MetricOrth}= \tilde{g}_{il}\,\tilde{g}_{jm}\,u^{ij}\,v^{lm}= \sum_{ij} h_{(i)}^2 h_{(j)}^2 u^{ij}v^{ij}$.

\subsubsection{Representation based on embedding}%
\label{sec:Embedding}
An alternative convenient representation, to be used in SFEM, TraceFEM, and DI, follows from considering the $n$-tensor fields as general mappings into the embedding space, e.g.,
\begin{equation}\label{eq:embedded-tensor-fields}
  {\ub}^{(0)} \colon \Surface \to \R,\quad {\ub}^{(1)} \colon \Surface \to \R^3,\quad {\ub}^{(2)} \colon \Surface \to L(\R^{3}, \R^{3})\,,
\end{equation}
with $L(\R^{3}, \R^{3})$ the linear mappings between $\R^3$ and $\R^3$, which can be represented as $\R^{3\times 3}$ tensor.

If we use the standard basis in $\R^3$, then the component vector $\TT{u}^{(n)}$, $n=1,2$, can be identified with the corresponding fields $\ub^{(n)}$. To simplify the notation, we use this identification and delete the underline in the notation of the tensor fields when the meaning is clear from the context.
This embedded representation leads to a natural inner product defined in the embedding, i.e., for $n$-tensors $\TT{u}^{(n)}\colonequals[u^{\alpha_1\ldots \alpha_n}]$ and $\TT{v}^{(n)}\colonequals[v^{\alpha_1\ldots \alpha_n}]$, we have
\[
  \scalarProd{\ub^{(n)}}{\vb^{(n)}} \colonequals u^{\alpha_1\ldots \alpha_n}\,v_{\alpha_1\ldots \alpha_n}\,,
\]
where the indices can be raised and lowered using the Euclidean metric $\delta_{\alpha\beta}$. Note that since the tensor fields are represented in the embedding space, the indices are in the range $\alpha_k\in\{1,2,3\}$.

Corresponding to the unit normal field $\nb$ we introduce the tangential projection $\P=\bm{I}-\nb\otimes\nb$ and denote in the following a general tensor projection operator for $n$-tensors $\TT{u}^{(n)}=[u^{\alpha_1\ldots \alpha_n}]$ as $\PP$, defined by componentwise projection,
\begin{equation}\label{eq:tensor-projection}
  [\PP\ub^{(n)}]\indices{^{\alpha_1\ldots \alpha_n}} \colonequals \tensor{P}{^{\alpha_1}_{\beta_1}}\cdots\tensor{P}{^{\alpha_n}_{\beta_n}} u\indices{^{\beta_1\ldots \beta_n}}\,.
\end{equation}

The (total) covariant derivative of tangential tensor fields $\ub^{(n)}$ with embedded representation $\Emb{\ub}^{(n)}$ can be defined as $\Grad_\Surface\ub^{(n)} \colonequals \PP\grad\Emb{\Ext{\ub}}^{(n)}$. Recall that the overline symbol denotes a (smooth) extension of a function to a surface neighborhood, while the underline symbol emphasizes that the Euclidean gradient $\nabla$ is applied componentwise. Written out for $n=0,1,2$, this definition reads:
\begin{align}
  \Grad_\Surface \ub^{(0)} &= \P\grad \Emb{\Ext{\ub}}^{(0)}, &&  \label{defGa}\\
  \Grad_\Surface \ub^{(1)} &= \P\grad \Emb{\Ext{\ub}}^{(1)}\P, && \label{defGb} \\
  \left[\Grad_\Surface\ub^{(2)}\right]\indices{^{\alpha_1 \alpha_2 \alpha_{3}}} &= \tensor{P}{^{\alpha_1}_{\beta_1}}\tensor{P}{^{\alpha_2}_{\beta_2}} \tensor{P}{^{\alpha_3}_{\beta_3}}\tensor{\delta}{^{\beta_3}^{\gamma}}\D{\Ext{u}^{\beta_1 \beta_2}}{x^{\gamma}}\,. && \label{defGc}
\end{align}
The definitions in \cref{defGa,defG1} give the same surface gradient operator for scalar functions. The definitions in \cref{defGb,defGc} are also used for vector and tensor fields that are not in the tangent bundle. In the case that these fields are tangential, these definitions yield the same gradient operators as those defined in \cref{defG2,defG3}.

The divergence of a tangential $n$-tensor, $n\geq 1$, is given by
\begin{equation} \label{defdiv}
 \left[ \Div_{\Surface} {\ub}^{(n)}\right]^{\alpha_1\ldots \alpha_{n-1}}
 =\tensor{P}{^{\alpha_1}_{\beta_1}}\cdots\tensor{P}{^{\alpha_{n}}_{\beta_{n}}}
 \frac{\partial \Ext{u}^{\beta_1 \ldots \beta_n}}{\partial x^{\alpha_{n}}}\,.
\end{equation}

\subsection{Problem definition}%
\label{sec:ProblemDefinition}
We study the following model problem: Find tangential $n$-tensor fields $\ub$ that solve
\begin{equation}\label{eq:strong-formulation}
  \partial_t\ub - \Laplace_\Surface\ub = 0\quad\text{ on }\Surface\,,
\end{equation}
subject to appropriate ``no-flux'' boundary conditions and initial conditions. Note that for $n=0$ the tangential condition is void. For $n \geq 1$, the $\Laplace_\Surface$ operator is the (negative) connection-Laplacian, the natural extension of the Laplace-Beltrami operator to $n$-tensor fields. It can be written as $\Laplace_\Surface\ub = \operatorname{div}_{\Surface}\Grad_\Surface\ub$, where $\Grad_\Surface$ is the covariant gradient operator defined above and $\Div_{\Surface}$ is the tensor surface divergence as in \cref{defdiv}.

For $n=0$ the PDE \eqref{eq:strong-formulation} corresponds to the scalar heat diffusion problem on a surface, which shares several properties with the corresponding equation in flat space. For example, it holds that $\Avg{\ub}(t) = \Avg{\ub^0}$ and $\ub(t,\xb) \to \Avg{\ub^0}$ for $t \to \infty$, with the mean $\Avg{\ub}(t) = \frac{1}{\text{area}(\Surface)} \int_{\Surface} \ub(t, \xb) \, d \xb$ of $\ub(t,\xb)$. It also holds that if $\ub^0 \geq 0$ on $\Surface$ and $\ub^0 \not\equiv 0$, then $\ub(t) > 0$ on $\Surface$ for $t > 0$, see for example \cite{PhDBuergerRegensburg}.
There are also results available that explain certain influences of surface curvature on the solution of the scalar heat equation problem. We outline some of these results.
Let $\ub^0(\xb) = \delta_{\pb}(\xb)$ for some $\pb\in\Surface$ be the Dirac delta function. In \cite{Varadhan1967Behavior} it is shown that for this initial condition the corresponding solution satisfies
\[
  \lim_{t \to 0} \big[-2 t \log \ub(t, \xb) \big] = d_{\Surface}^{\,2}(\xb,\pb)\,,
\]
with $d_{\Surface}(\cdot, \cdot)$ the geodesic distance on $\Surface$. This property is  used to approximate geodesic distances on curved surfaces in computer science \cite{CW2013Geodesics}. In \cite{Faraudo2002Diffusion} the scalar heat equation is considered and rewritten in terms of geodesic polar coordinates. This allows to separate the diffusion from the geometric influence, the latter being completely determined by the Gaussian curvature $\GaussCurvature$. We will now outline a result that will be used to explain a certain phenomenon observed in the numerical experiments in \Cref{sec:NumericalExperiments}. The solution of the scalar heat equation with initial condition $\delta_{\pb}$ can be expressed as
\begin{equation}\label{eq:kernel-folding}
  \ub(t, \xb) = \int_{\Surface} k_t(\xb,\yb) \delta_{\pb}(\yb) \; d \yb
\end{equation}
with heat kernel $k_t(\cdot,\cdot)$. In \cite{SOG2009Concise,MS1967Curvature} it is shown that
\begin{equation}\label{eq:heat-kernel}
  k_t(\xb,\xb) = \frac{1}{4 \pi t} \big(1 + \frac{1}{6} \GaussCurvature(\xb)\, t + \text{h.o.}(t) \big)\,.
\end{equation}
This result motivates general statements such as ``heat tends to diffuse slower at points of positive curvature and faster at points of negative curvature''. Several analytical solutions of the heat equation for special surfaces also have been derived \cite{Faraudo2002Diffusion}.

For $n = 1$, the surface vector heat equation, some of these results can be generalized. In particular, it can be shown, again by considering the associated heat kernel, that for $t \to 0$ it behaves like parallel transport along geodesics, along with a decay in magnitude that is identical to the decay of the scalar heat kernel \cite{SS2019Vector}. This property is crucial for several applications in computer graphics \cite{KCPS2013Globally,SS2019Vector} and data science \cite{SW2012Vector}, where it is also extended to tensor fields with $n > 1$. Since the tangential tensor-valued heat equation can also be seen as the $L^2$-gradient flow for the tensor Dirichlet energy $\int_{\Surface} \Norm{\Grad_\Surface \ub }^2 \,\textrm{d}\xb$, its solution tends to the minimizer of this energy functional (``smoothest possible'' tensor-field)  for $t\to\infty$.

\section{Finite element discretization schemes}%
\label{sec:FiniteElementDiscretizationSchemes}
To be able to  apply a finite element discretization method we consider the $n$-tensor diffusion problem in a variational setting, using standard notation for Bochner spaces:

\begin{problem}\label{prob:Problem1}
  Find $\ub\in C^1(0,T;\bm{H}^1(\Surface, \TensorBundle^n\Surface))$ such that
  \begin{equation}\label{eq:variational-formulation}
    \Inner{\partial_t\ub(t)}{\vb}_{\Surface} + \Inner{\Grad_{\Surface}\ub(t)}{\Grad_{\Surface}\vb}_{\Surface} = 0\,,\quad\text{for all}~ \vb\in\bm{H}^1(\Surface, \TensorBundle^n\Surface)\,,
  \end{equation}
  for $t\in(0,T]$, subject to $\ub(0) = \ub^0$. Here $\Inner{\cdot}{\cdot}_{\Surface}$ denotes the (tensor) $L^2$-scalar product.
\end{problem}

The ISFEM is based directly on the variational formulation given in \Cref{prob:Problem1}. The other three methods, SFEM, TraceFEM, and DI, use the surface embedding in $\R^3$ and gradient representations as presented in \Cref{sec:Embedding}. For these methods applied to the $n$-tensor problem with $n\geq 1$, it is natural to allow (small) \emph{non}tangential solution components. The variational formulation given in \Cref{prob:Problem1} is not a suitable starting point for such a finite element method, since it uses the range space $\TensorBundle^n\Surface$.
We now introduce, for $n \geq 1$, an augmented variational formulation with range space $\TensorBundle^n\R^3 \simeq\R^{3^n}$. It uses a term $\Inner{\QQ\ub(t)}{\QQ\vb}_{\Surface}$,
where $\QQ=\textrm{Id}-\PP$ is the normal projection operator that is scaled by  a penalty parameter $\omega > 0$. Note that this term vanishes for tangential functions. The augmented variational formulation reads as follows:

\begin{problem}\label{prob:Problem2}
  Assume $n \geq 1$. Find $\ub\in C^1(0, T; \bm{H}^1(\Surface, \TensorBundle^n\R^3))$ such that
  \begin{equation}\label{eq:variational-formulation-embedded}
    \Inner{\partial_t\PP\ub(t)}{\PP\vb}_{\Surface}
    + \Inner{\Grad_{\Surface}\PP\ub(t)}{\Grad_{\Surface}\PP\vb}_{\Surface}
    + \omega\,\Inner{\QQ\ub(t)}{\QQ\vb}_{\Surface} = 0
  \end{equation}
  for all $\vb\in\bm{H}^1(\Surface, \TensorBundle^n\R^3)$ and $t\in(0,T]$, with initial condition $\ub(0)=\ub^0$.
\end{problem}

Note that in \cref{eq:variational-formulation-embedded} first derivatives appear only for the tangential components $\PP\ub$, but not for the normal components $\QQ\ub$.
\Cref{prob:Problem2} is consistent with \Cref{prob:Problem1} in the following sense. Let $\ub_1$ and $\ub_2$ be two solutions of \Cref{prob:Problem2} and $\bm{w}\colonequals\ub_1-\ub_2$. Then we have $\bm{w}(0)=0$ and from \cref{eq:variational-formulation-embedded} we obtain $\partial_t\|\PP\bm{w}\|_{\Surface}^2 \leq 0$ for $t \in [0,T]$. Hence $\PP\bm{w}(t)=0$ for $t \in [0,T]$. Using this in \cref{eq:variational-formulation-embedded} it follows that $\QQ\bm{w}(t)=0$ and thus $\bm{w}(t)=0$ for $t \in [0,T]$. We conclude that we have a unique solution of \Cref{prob:Problem2}. It is easy to verify that a solution of \Cref{prob:Problem1} is also a solution of \Cref{prob:Problem2}.  We see that by adding the consistent penalty term, with $\omega \geq 0$ arbitrarily, we do not change the continuous solution.  In general, this ``exact'' consistency property does not hold after discretization, and one must then choose an appropriate value of the penalty parameter to control the consistency error, see Sections \ref{sec:sfem} and \ref{sec:tracefem}.

In SFEM and TraceFEM, which will be introduced below, a discrete projection operator $\PP_h$ is used. This projection operator is generally discontinuous across element boundaries. Thus, for a vector- or tensor-valued finite element function $\ub_h$, the projected function $\PP_h\ub_h$ has no global $H^1$-smoothness, and applying a discrete gradient to it is a nontrivial problem. To circumvent this issue, we apply the product rule to  the term  $\Grad_{\Surface}\PP\ub$ in \cref{eq:variational-formulation-embedded} as follows. We have $\Grad_{\Surface}\PP\ub=\Grad_{\Surface}\ub- \Grad_{\Surface}\QQ\ub$. For a vector field $\ub= \ub^{(1)}$ it is easy to check that $\Grad_{\Surface}\QQ\ub =\Grad_{\Surface}( (\nb\otimes\nb) \ub )= - \scalarProd{\ub}{\nb}\bm{H}$ holds. So we get $\Grad_{\Surface}\PP\ub=\Grad_{\Surface}\ub +\scalarProd{\ub}{\nb}\bm{H}$ and the representation on the right-hand side is suitable for a finite element approximation. We introduce the notation $G(\ub)\colonequals\scalarProd{\ub}{\nb}\bm{H}$. Note that $G$ depends on the extended Weingarten map. Similar results hold for $n \geq 2$.
In tensor notation we get the following identities:
\begin{align}\label{eq:gradQu}
  [\Grad_{\Surface}\QQ {{\ub}}^{(1)}]\indices{^{\alpha_1 \alpha_2}}= [G(\ub)]\indices{^{\alpha_1 \alpha_2}}
    &\colonequals -H\indices{^{\alpha_1 \alpha_2}}\,u^\beta\,n\indices{_\beta}, && n=1\,, \\
  [\Grad_{\Surface}\QQ {{\ub}}^{(2)}]\indices{^{\alpha_1 \alpha_2 \alpha_3}}=[G(\ub) ]\indices{^{\alpha_1 \alpha_2 \alpha_3}}
    &\colonequals -H\indices{^{\alpha_1 \alpha_3}}\,P\indices{^{\alpha_2}_{\beta}}u\indices{^{\gamma_1 \beta}}\,n\indices{_{\gamma_1}}
        -H\indices{^{\alpha_2 \alpha_3}}\,P\indices{^{\alpha_1}_{\beta}}u\indices{^{\beta \gamma_2}}\,n\indices{_{\gamma_2}}, && n=2\,. \notag
\end{align}
Thus we have the following alternative representation of the second term in \cref{eq:variational-formulation-embedded}, which will be the one used in SFEM and TraceFEM below:
\begin{align}\label{formeq}
  \Inner{\Grad_{\Surface}\PP\ub(t)}{\Grad_{\Surface}\PP\vb}_{\Surface}
    &= \Inner{\Grad_{\Surface}\ub(t) + G(\ub(t)) }{\Grad_{\Surface}\vb + G(\vb)}_{\Surface}\,.
\end{align}

In the following subsections, we briefly discuss four well-known finite element discretization methods and apply them to the spatial discretization of the $n$-tensor heat problem. We combine these spatial discretizations with a standard BDF-2 time discretization. We start with  ISFEM, which is the ``most conforming'' method in the sense that it is based on the variational formulation in \Cref{prob:Problem1}. This method uses intrinsic gradient representations. The methods SFEM, TraceFEM, and DI are based on the formulation in \Cref{prob:Problem2} and use gradient representations in the embedding space. The DI method considers an additional approximation of the inner products by ``extending'' the PDE in a domain $\Omega \subset \R^3$ containing the surface $\Surface$ and numerically restricting the integrals to $\Surface$ using a smeared-out Dirac-delta function. This method is not consistent in the sense that the solution of the extended PDE, restricted to $\Surface$, does \emph{not} coincide with the solution of Problems 1 and 2. In this sense, the DI approach is the ``least conforming'' one. In the presentation of the methods we restrict ourselves to the case of the lowest order finite elements. In remarks we will briefly comment on extensions to higher order finite elements.

In \Cref{sec:Comparison} we discuss and compare the four methods and in particular address the issues a)--c) formulated in the introduction.

\subsection{Intrinsic Surface Finite Element Method (ISFEM)}%
\label{sec:isfem}
The ISFEM has been introduced only for scalar-valued problems in \cite{BFP2021Intrinsic}. We briefly review the scalar ISFEM setting and extend it to the case of a vector-valued problem. An analogous extension to tensor fields is possible, but it has not yet been addressed. The main idea of ISFEM is to consider the formulation as given in \Cref{prob:Problem1} and to discretize it in local coordinates with the intrinsic differential operators, which contain all geometric information. We will consider the local coordinates $\sv=(\hat{s}^1,\hat{s}^2)$ with respect to the orthogonal tangent vectors $\BasisOrth_{1},\BasisOrth_{2}$, and the intrinsic differential operators defined using \cref{eq:crhistsymb-orth}.

Let  $\Tri{\Surface}$ be a (curved) exact surface triangulation, formed by a set of non-intersecting (curved) surface triangles with vertices on $\Surface$, such that $\Surface = \bigcup_{\Elem\in\Tri{\Surface}}\Elem$.
We will introduce \emph{conforming} subspaces $\VVh^{(n)}_{\Surface}$, so that the relation $\VVh^{(n)}_{\Surface} \subset \bm{H}^1(\Surface, \TensorBundle^n\Surface)$ holds, see \Cref{prob:Problem1}. These spaces are used in an approximate Galerkin discretization of \Cref{prob:Problem1}, in the sense that the surface integrals $\Inner{\cdot}{\cdot}_\Elem$ are approximated by a quadrature rule.

We  denote by $\Inner{\ub}{\vb}_{h}\colonequals\sum_q w_q \scalarProd{\ub(\xb_q)}{\vb(\xb_q)}_{\MetricOrth(\xb_q)} \sqrt{\Abs{\MetricOrth(\xb_q)}}$ such an approximation of the surface integrals $\Inner{\ub}{\vb}_{\Elem}$ by a quadrature rule with $\xb_q\in\Elem$ the quadrature points and $w_q\in\R$ the associated quadrature weights.
In terms of practical computation, the key point is the need for geometric information only at quadrature points, in an exact or approximate way.
\begin{remark}
  In the benchmark problem considered, see \Cref{sec:NumericalExperiments}, we apply a Gauss quadrature rule of order three. In this case, we use the knowledge of the surface parametrization to assign geometric information at the quadrature points.
\end{remark}

First we consider the scalar case $n=0$. The function space $\VVh^{(0)}_{\Surface}=\operatorname{span}(\{\psi_{l}\})$ is spanned by continuous basis functions $\psi_l\colon\Surface\to\R$ that are obtained by formally gluing together localized functions $\psi_l^{\Elem}\colon\Elem\to\R$ for $\Elem\in\Tri{\Surface}$.
For each element $\Elem\in\Tri{\Surface}$, we consider the associated element $\hat{\Elem}=\ChartMap^{-1}(\Elem)$ in the reference domain and we define the classical linear Lagrange nodal basis functions $\hat{\psi}_l^{\hat{\Elem}}(\xv)$ in reference local coordinates  $\xv\in\hat{\Elem}$. If we denote by $\xb=\ChartMap(\xv)\in\Elem$ the corresponding associated surface coordinates in $\Elem$, then the surface basis functions are simply lifted using this mapping, i.e., $\psi_l^{\Elem}(\xb) \colonequals \hat{\psi}_l^{\hat{\Elem}}(\xv)$.
In order to compute gradients in the tangential basis representation $\{\BasisOrth_{1},\BasisOrth_{2}\}$ associated with coordinates $\sv$ instead of the natural tangential basis $\{\Basis_1,\Basis_2\}$ associated with $\xv$, we need to perform a coordinate transformation, i.e.,
\begin{equation*}
  \grad_{\Surface}\psi_l^{\Elem}(\xb) \colonequals \MetricOrth^{-1}\bm{W}\hat{\grad}\hat{\psi}_l^{\hat{\Elem}}(\xv)\,,
\end{equation*}
where $\bm{W}=\tilde{\bm{J}}^{+}\bm{J}$ is the Jacobian of the coordinate change between $\xv$ and $\sv$, with $\tilde{\bm{J}}=[\BasisOrth_{1},\BasisOrth_{2}]$, $\tilde{\bm{J}}^{+}$ its pseudoinverse, and $\bm{J}=[\Basis_1,\Basis_2]$.
\begin{remark}
  In the case of a surface obtained by the graph of a scalar function, for example $\ChartMap(\xv)=(\hat{x}^1,\hat{x}^2,f(\hat{x}^1,\hat{x}^2))^T=\xb\in\Surface$, the matrix $\bm{W}$ is obtained directly from the $2\times 2$ block of $\tilde{\bm{J}}^T$ corresponding to the independent variables.
\end{remark}

The discrete scalar functions $\ub_h^{(0)}\in\VVh^{(0)}_{\Surface}$ can be expanded in terms of the basis functions as $\ub_h^{(0)}(\xb)=\sum_l u_l \psi_{l}(\xb)$, where $u_l$ is the scalar coefficient associated with the basis function $\psi_{l}$. For discrete vector-valued functions $\ub_h^{(1)}$ we use the orthogonal covariant reference frame $\{\BasisOrth_{1}{},\BasisOrth_{2}{}\}$ and represent the solution in contravariant components, i.e., $\Emb{\ub}^{(1)}=\left[u^{i}\right]$, for $i=1,2$.
Each component $u^i$ can be approximated by the discrete functions $u_h^i=\sum_l u_l^i\,\psi_l$, where $\{\psi_l\}_l$ is the set of scalar basis functions of $\VVh^{(0)}_{\Surface}$.
Thus we get
\begin{equation*}
  \ub^{(1)}_h(t) = \sum_l u_l^1(t)\,\psi_l \BasisOrth_{1} + u_l^2(t)\,\psi_l\BasisOrth_{2}\,,
\end{equation*}
which gives rise to the definition of a discrete vector function space:
\begin{equation*}
  \VVh^{(1)}_{\Surface} \colonequals \Set{\vb_h=v_h^1\,\BasisOrth_{1} + v_h^2\,\BasisOrth_{2} \mid v_h^1,v_h^2\in\VVh^{(0)}_{\Surface}}\,.
\end{equation*}
The same idea can be used to define a discrete tensor function space $\VVh^{(n)}_{\Surface}$.

By applying the definition of gradients and scalar product in~\Cref{secintri}, with respect to the orthogonal reference frame $\{\BasisOrth_{1},\BasisOrth_{2}\}$, and the quadrature rule $\big(\cdot,\cdot\big)_{h}$, we obtain the semi-discrete ISFEM discretization of \Cref{prob:Problem1}:

\begin{problem}
  Find $\ub_h(t)=\ub^{(n)}_h(t)\in\VVh^{(n)}_{\Surface}$ such that
  \begin{equation}\label{eq:isfem-semi-discrete}
    \Inner{\partial_t\ub_h(t)}{\vb_h}_{h} + \Inner{\Grad_{\Surface}\ub_h(t)}{\Grad_{\Surface}\vb_h}_{h} = 0\,,\quad\text{for all}~ \vb_h\in\VVh^{(n)}_{\Surface}
  \end{equation}
  for $t\in(0,T]$, with initial condition $\ub_h(0)=\ub^0$.
\end{problem}

We obtain fully discrete schemes by applying a BDF-2 discretization scheme to the semi-discrete problem in \cref{eq:isfem-semi-discrete}.

\subsection{Surface Finite Element Method (SFEM)}%
\label{sec:sfem}
The essence of the lowest order SFEM is the approximation of  $\Surface$ by a (shape regular) triangulation, consisting of flat triangles, and the use of globally continuous piecewise linears on this triangulation to approximate the continuous solution. This technique avoids surface parametrizations and is very similar (for scalar problems) to a standard finite element method in a flat domain. The piecewise triangular surface approximation is denoted by $\Surface_h$.
The space of globally continuous piecewise linears on $\Surface_h$ is denoted by $V_{\Surface_h}$. On each triangle $\Elem$ of $\Surface_h$ we introduce the natural geometry normal $\nb_h\vert_{\Elem}$. These local normals $\nb_h\vert_{\Elem}$ are formally glued together to form the discrete surface normal field $\nb_h$. The discrete tangential projection is given by $\Ph=\bm{I}-\nb_h\otimes\nb_h$. The SFEM for the scalar case is well-known in the literature and reads as follows (cf. \Cref{prob:Problem2}), with $\Grad_{\Surface_h}\vb_h = \P_h\grad \Ext{\vb}_h$ the discrete analog of the surface gradient as in \cref{defGa}:

\begin{problem}[Scalar problem]\label{prob:discrSFEM1}
  Find $\ub_h=\ub_h^{(0)}\in C^1(0,T; V_{\Surface_h})$ such that
  \begin{equation}
  (\partial_t\ub_h(t), \vb_h )_{\Surface_h} + (\Grad_{\Surface_h}\ub_h(t), \Grad_{\Surface_h}\vb_h )_{\Surface_h} = 0\,,\quad\text{for all}~ \vb_h\in V_{\Surface_h}
  \end{equation}
  and for all $t\in(0,T]$ subject to an initial condition $\ub_h(0)=\bm{I}_h\ub^0$. Here $\bm{I}_h$ denotes the nodal interpolation operator in the finite element space $V_{\Surface_h}$.
\end{problem}

Using the nodal finite element basis in the space $V_{\Surface_h}$ results in an ODE system for the  coefficients of $\ub_h$.

We now consider $n \geq 1$. The discretization is based on the formulation in \Cref{prob:Problem2}, combined with a componentwise approximation using SFEM for scalar-valued problems. For a detailed description for tensor-valued problems see \cite{NN2019Finite,HL2020Analysis,HaPr2021Tangential}.

The discrete tensor projection operator $\PPh$ is defined analogously to \cref{eq:tensor-projection}. A corresponding orthogonal projection $\QQh=\mathcal{I}-\PPh$ follows naturally.
We use the surface finite element space $\VVh^{(n)}_{\Surface_h}=[\Vh_{\Surface_h}]^N$ as the product space of $N=3^n$ scalar Lagrange spaces.
A discrete surface gradient $\Grad_{\Surface_h}$ is defined as in \cref{defGb,defGc}, but with the continuous projections replaced by the discrete analogons. Error analysis and numerical experiments show  that replacing the projection operator $\QQ$ in the penalty term of the continuous variational formulation by its discrete analog $\QQ_h$ is \emph{not} satisfactory, since it leads to suboptimal convergence in the $L^2$-norm, cf. \cite{HL2020Analysis}. Optimal convergence is obtained by using instead a projection operator based on a normal $\nb^\sharp_h$, which is a \emph{one order more accurate approximation} of $\nb$ than the $\Surface_h$-normal $\nb_h$. We denote such a modified (``higher order'') projection by $\QQ_h^\sharp$.

Thus we obtain the following SFEM discretization of \cref{eq:variational-formulation-embedded}, where we use the result \eqref{formeq}, see also \cite{HaPr2021Tangential}:

\begin{problem} \label{prob:SFEMdiscr2}
  Take $n \geq 1$. Find $\ub_h=\ub_h^{(n)} \in C^1(0,T;\VVh^{(n)}_{\Surface_h})$ such that
  \begin{multline}\label{eq:sfem-discretization}
    \big(\partial_t\PPh {\ub}_h(t), \PPh {\vb}_h\big)_{\Surface_h} + \big(\Grad_{\Surface_h} \ub_h(t) +G_h(\ub_h (t)), \Grad_{\Surface_h}{\vb}_h + G_h(\vb_h)\big)_{\Surface_h} \\
    + \beta h^{-2}\,\big(\QQ_h^\sharp {\ub}_h(t), \QQ_h^\sharp {\vb}_h\big)_{\Surface_h} = 0 \quad \text{for all}~{\vb}_h\in\VVh^{(n)}_{\Surface_h}
  \end{multline}
  and for all $t\in(0,T]$ subject to an initial condition ${\ub}_h(0)=\bm{I}_h {\ub}^0$.
\end{problem}

Here $G_h(\cdot)$ is a discrete analog of $G(\cdot)$ in \cref{formeq}, e.g., for $n=1$, $G_h(\vb_h)=\scalarProd{\vb_h}{\nb_h}\bm{H}_h$, where $\bm{H}_h$ is an approximation of the Weingarten mapping. The parameter $\beta > 0$ is a penalty parameter. The scaling with $h^{-2}$ in the penalty term follows from an error analysis, cf. \cite{HL2020Analysis,HaPr2021Tangential}.

\begin{remark}
  The discrete Weingarten map $\bm{H}_h$ can be computed from the elementwise gradient of the discrete normal field $\nb_h$. Using the representation $\nb_h=\mb_h/\|\mb_h\|$, with $\mb_h$ being the cross-product of the columns of $\bm{J}_h$ and thus a discrete function, we can compute $\bm{H}_h=\PPh\grad\nb_h=\|\mb_h\|^{-1}\PPh\grad (I_h \mb_h)$.
\end{remark}

\begin{remark}
  If the surface $\Surface$ is described by the coordinate mapping $\ChartMap$, higher order surface approximations than piecewise linear can be obtained by (Lagrange) interpolation $\ChartMap_h=I_h\ChartMap$. With discrete functions defined in the reference domain and lifted to the discrete surface using $\ChartMap_h$, a higher-order function space $\Vh_{\Surface}$ can be constructed, cf. \cite{Demlow2009HigherOrder}. Similar to the piecewise flat surface and linear function setting, geometric quantities are obtained by derivatives of the discrete parametrization $\ChartMap_h$. This also allows for high-order convergence for the projection based scheme, cf. \cite{HL2020Analysis,HaPr2021Tangential}.

  If derivatives of the continuous parametrization $\ChartMap$ are directly available and computable, an exact parametrization of the surface geometry is also possible, cf. \Cref{sec:isfem}, and is used in the numerical example to compute a reference solution.
\end{remark}

The discretization in time follows standard approaches and is therefore not described in detail. We consider a classical BDF-2 scheme.

\subsection{Trace Finite Element Method (TraceFEM)}%
\label{sec:tracefem}
The TraceFEM is based on the same variational \Cref{prob:Problem2} as the SFEM.  The former  uses a finite element space which is defined on a background volumetric mesh  that is not fitted to the surface. The geometry approximation is based on an implicit description of the surface using a level set approach. For an overview of TraceFEM we refer to \cite{OlshanskiiReusken2017}.

We assume that the surface $\Surface$ is represented as the zero level of a level set function $\phi$. We denote by $\Domain$ a sufficiently small polygonal 3d neighborhood of the surface. The surface approximation is based on a piecewise linear approximation $\phi_h$ (e.g., linear interpolation) of $\phi$ and is given by $\Surface_h\colonequals\Set{\xb \in \Omega \mid \phi_h(\xb)=0}$. Let $\Tri{\Domain}$ be a shape regular tetrahedral triangulation of $\Domain$ and $\Vh_{\Domain}$ be the standard finite element spaces of continuous piecewise linear polynomials  on $\Tri{\Domain}$.  For higher order constructions of $\Surface_h$ see \Cref{rem_TraceFEM_surface_approximation}.
We introduce the set $\Tri{\Domain}^{\Surface_h}$, which consists of all tetrahedra $\Elem\in\Tri{\Domain}$ that have a nonzero intersection with $\Surface_h$. The domain formed by all these tetrahedra is denoted by $\Domain^{\Surface_h}=\bigcup_{\Elem\in\Tri{\Domain}^{\Surface_h}}\Elem$. On $\Domain^{\Surface_h}$ we define by simple restriction the scalar finite element space $\Vh^{\Surface_h}_{\Domain}\colonequals \Set{v\vert_{\Domain^{\Surface_h}} \mid v\in \Vh_{\Domain}}$. A corresponding $n$-tensor finite element space is given by $\VVh^{\Surface_h}_{\Domain} \colonequals \big[ \Vh^{\Surface_h}_{\Domain} \big]^N$ with $N=3^n$.
To avoid instabilities due to small cuts of $\Surface_h$ in the triangulation $\Domain^{\Surface_h}$, a so-called normal derivative volume stabilization is used \cite{BHLM2018Cut, GLR2018Analysis}. Again, $\nb_h$ denotes the piecewise normal field on $\Surface_h$, $\bm{P}_h = \bm{I} - \nb_h \otimes \nb_h$, and the discrete tensor projection operator $\PP_h$ is defined as in Subsection \ref{sec:sfem}.

We now describe the method for the scalar case $n=0$. The stabilization is then given by $s_h(\ub_h,\vb_h)\colonequals\big(\nb_h \cdot \nabla \ub_h, \nb_h \cdot \nabla \vb_h)_{\Domain^{\Surface_h}}$. The discrete problem is as follows, see \Cref{prob:discrSFEM1}:

\begin{problem}[Scalar problem]\label{prob:discrTraceFEM1}
  Find $\ub_h=\ub_h^{(0)}\in C^1(0,T; \Vh^{\Surface_h}_{\Domain})$ such that
  \begin{equation}
  \big(\partial_t\ub_h(t), \vb_h \big)_{\Surface_h} + \big(\Grad_{\Surface_h}\ub_h(t), \Grad_{\Surface_h}\vb_h \big)_{\Surface_h} + \beta' h^{-1} s_h(\ub_h(t),\vb_h) = 0
  \end{equation}
  for all $\vb_h\in  \Vh^{\Surface_h}_{\Domain}$ and all $t\in(0,T]$, subject to an initial condition ${\ub}_h(0)=\textbf{I}_\Domain^{\Surface_h}\Ext{\ub}^0$. Here $\textbf{I}_\Domain^{\Surface_h}$ denotes the nodal interpolation operator in the finite element space $\Vh^{\Surface_h}_{\Domain}$.
\end{problem}

Note that compared to \Cref{prob:discrSFEM1} we use a different finite element space and have added the stabilization term $s_h(\cdot,\cdot)$.
This stabilization term significantly improves the conditioning of the resulting linear systems. In case of a smooth closed surface, the condition number of the stiffness matrix corresponding to the Laplace-Beltrami operator has the usual $h^{-2}$ growth. For smooth surfaces with boundary, very strong ill-conditioning can still occur in certain situations. We will not discuss this effect, which does not occur for the surface considered in \Cref{sec:NumericalExperiments}.

We now consider $n \geq 1$. In the same spirit as in the SFEM, cf. \Cref{prob:SFEMdiscr2}, it is based on the variational formulation in \Cref{prob:Problem2}. The normal derivative volume stabilization has the form
\begin{equation*}
  \bm{s}_h(\ub,\vb) \colonequals \int_{\Domain^{\Surface_h}} \left(\grad\ub \Tdot{n+1} \nb_h \right) \cdot \left(\grad\vb \Tdot{n+1} \nb_h \right) \,\text{d}\xb\,,
\end{equation*}
with $[\grad\ub \Tdot{n+1} \nb]\indices{^{\alpha_1,\dots,\alpha_n}}=\frac{\partial\bar{u}\indices{^{\alpha_1,\ldots,\alpha_n}}}{\partial x\indices{^{\alpha_{n+1}}}}n\indices{^{\alpha_{n+1}}}$ and we obtain the following discretization of \cref{eq:variational-formulation-embedded}:

\begin{problem} \label{prob:TraceFEMdiscr2}
  Take $n\geq 1$.	Find ${\ub}_h={\ub}^{(n)}_h \in C^1(0,T;\VVh^{\Surface_h}_{\Domain})$ such that
  \begin{equation} \label{eq:tracefem-time-discretization}
    \begin{split}
    & \big(\partial_t\PPh {\ub}_h(t), \PPh {\vb}_h\big)_{\Surface_h}
    + \big(\Grad_{\Surface_h} \ub_h(t) +G_h(\ub_h (t)), \Grad_{\Surface_h}{\vb}_h+ G_h(\vb_h)\big)_{\Surface_h}
    \\ & \quad + \beta h^{-2}\,\big(\QQ_h^\sharp {\ub}_h(t), \QQ_h^\sharp {\vb}_h\big)_{\Surface_h} + \beta' h^{-1}\, \bm{s}_h(\ub_h(t),\vb_h) = 0 \quad \text{for all}~{\vb}_h\in\VVh^{\Surface_h}_{\Domain}
    \end{split}
  \end{equation}
  and for all $t\in(0,T]$ subject to an initial condition ${\ub}_h(0)={\textbf{I}}_\Domain^{\Surface_h} \Ext{\ub}^0$, with $\beta > 0$ a penalty parameter and $\beta' > 0$ a stabilization parameter.
\end{problem}

As in \Cref{sec:sfem} we use an ``improved'' projection $\QQ_h^\sharp$ based on a higher  order normal approximation. A motivation  for this improved projection and a construction of an improved normal approximation are given in \cite{JankuhnReusken2020}. As in SFEM the term $G_h(\cdot)$ is a discrete analog of $G(\cdot)$ given in \cref{formeq}. The semi-discrete \Cref{prob:TraceFEMdiscr2} is essentially the same as the SFEM \Cref{prob:SFEMdiscr2} except for the additional stabilization term $\bm{s}_h(\cdot,\cdot)$.

As in \Cref{sec:sfem}, we use a classical BDF-2 scheme for time discretization.

\begin{remark}\label{rem_TraceFEM_surface_approximation}
  To obtain a higher order discretization method an isoparametric mapping $\Theta_h$ is the key ingredient. The main idea and construction of this mapping is explained in \cite{LehrenfeldReusken2017}. It is based on a level set function approximation $\phi_h \in \Vh^k_{\Domain}$ of order $k$. This function implicitly defines a surface approximation. For $k \geq 2$ numerical integration is difficult to implement. To obtain a computationally efficient method a piecewise triangular surface approximation $\Surface^{\mathrm{lin}}$ is used, defined as follows. Let $\hat{\phi}_h = I^1 \phi_h$ be the linear nodal interpolation  of the higher order level set function approximation $\phi_h$. Based on this, we define
  \begin{equation*}
    \Surface_{h} \colonequals \Theta_h(\Surface^{\mathrm{lin}}) = \Set{\xb \mid \hat{\phi}_h(\Theta_{h}^{-1}(\xb)) = 0}.
  \end{equation*}
  In the same way, the parametric mapping induces (higher order) finite element spaces.
\end{remark}

\subsection{Diffuse-Interface Approach (DI)}
The DI method, see \cite{RV2006PDEs, LLTVW2009DiffuseInterface, NNPV2018Orientational}, considers an approximation of \cref{eq:variational-formulation-embedded}, which is a classical problem in the embedding space $\R^3$ and thus leads to a setup where established standard volume FEM can be applied. Similar to \Cref{sec:tracefem}, the geometry approximation is based on an implicit description of the surface, but instead of a level set approach a phase field description is used. We define
\[
  \phi_\epsilon(\xb) \colonequals \frac{1}{2}\left(1 - \tanh\left(\frac{3}{\epsilon}\rho(\xb)\right)\right),\quad
  \delta_\epsilon(\xb) \colonequals \frac{36}{\epsilon}\phi_\epsilon^2(\xb)(1 - \phi_\epsilon(\xb))^2\,,
\]
for $\xb\in U_\delta(\Surface)$ with $0<\epsilon<\delta$ an interface thickness parameter. The phase-field function $\phi_\epsilon$ is based on the signed-distance representation $\rho$ of $\Surface$. With this definition of $\phi_\epsilon$ we get $\delta_\epsilon \to \delta_{\Surface}$ for $\epsilon \to 0$, with $\delta_{\Surface}$ the surface delta-function to $\Surface$. We define an extension of scalar-valued fields defined on $\Surface$ to the neighborhood $U_\delta(\Surface)$ by using the closest point projection, $\Ext{f}(\xb) \colonequals f(\pi(\xb)) = f(\xb - \rho(\xb) \nabla \rho(\xb))$ for $\xb\in U_\delta(\Surface)$. In the rest of the domain $\Omega$, the function $f$ is extended in an approximate way, e.g., using fast-marching algorithms or a Hopf--Lax algorithm \cite{ClaudelBayen2010LaxHopf}. The scalar-valued phase-field function $\phi_\epsilon$ and the delta-function $\delta_\epsilon$ are extended with a constant value.
Vector and tensor fields are extended by a componentwise extension of the embedded description.

As in \Cref{sec:tracefem} let $\Tri{\Domain}$ be a shape regular tetrahedral triangulation of $\Domain$ and $\Vh_{\Domain}$ be the standard finite element spaces of continuous piecewise linears defined on $\Tri{\Domain}$. We define the tensor finite element space $\VVh^{(n)}_{\Domain} \colonequals \big[\Vh_{\Domain}\big]^N$ as the product of $N=3^n$ scalar finite element spaces.

For scalar fields $\ub^{(0)}$, the FEM discretization of the DI approximation of \cref{eq:variational-formulation} reads:

\begin{problem}[Scalar diffuse interface approach \cite{RV2006PDEs}]
  Find $\ub_h = \ub_h^{(0)} \in C^1(0,T; \Vh_{\Domain})$ such that
  \begin{equation}\label{eq:di-scalar-time-discretization}
    \sum_{\Elem\in\Tri{\Domain}}\Inner{\delta_\epsilon \partial_t \ub_h(t)}{\vb_h}_{\Elem} + \Inner{\delta_\epsilon \grad \ub_h(t)}{\grad\vb_h}_{\Elem} + \sigma \Inner{ \left(1 - C_\epsilon \delta_\epsilon \right) \grad \ub_h(t)}{\grad
  \vb_h}_{\Elem} = 0
  \end{equation}
  for all $\vb_h \in\Vh_{\Domain}$ and for all $t\in(0,T]$, subject to the initial condition $\ub_h(0)=\Ext{\ub}^0$. As volume stabilization we add a small amount of ``additional diffusion'', with $\sigma = 10^{-8}$ in the numerical examples. The domain off the interface is characterized by $\left(1 - C_\epsilon \delta_\epsilon \right)$, with $C_\epsilon = 1/\max(\delta_{\epsilon})=4/9\,\epsilon$.
\end{problem}

For vector fields $\ub^{(1)}$ we consider the componentwise reformulation of the surface problem as in \cref{eq:variational-formulation-embedded}. For such a formulation we apply the scalar DI approach for each component. This requires geometric properties of the surface, namely the normal $\nb$ and the curvature $\bm{H}$ in the $\epsilon$-neighborhood of $\Surface$. To obtain these quantities one can use a numerical approximation of the signed-distance function $\rho^\epsilon$ in $\Omega$ and define $\nb^{\epsilon} \colonequals \grad \rho^\epsilon$ and $\bm{H}^{\epsilon} \colonequals -\nabla^2 \rho^\epsilon$, as well as the corresponding projections $\bm{P}^{\epsilon}$, $\PP^{\epsilon}$, and $\QQ^{\epsilon}$ with respect to $\nb^\epsilon$.

\begin{problem}
  Assume $n\geq 1$. Find $\ub_h\in C^1(0,T;\VVh^{(n)}_{\Domain})$ such that
  \begin{multline}\label{eq:di-discretization}
    \sum_{\Elem\in\Tri{\Domain}}\Inner{\delta_\epsilon \partial_t\PP^{\epsilon}\ub_h(t)}{\PP^{\epsilon}\vb_h}_{\Elem}
    + \Inner{\delta_\epsilon \Grad_{\Elem}\PP^{\epsilon}\ub_h(t)}{\Grad_{\Elem}\PP^{\epsilon}\vb_h}_{\Elem}
    + \beta\,\Inner{\delta_\epsilon\QQ^{\epsilon}\ub_h(t)}{\QQ^{\epsilon}\vb_h}_{\Elem} \\
    +\sum_{\Elem\in\Tri{\Domain}}\sigma\,\Inner{\left(1 - C_\epsilon \delta_\epsilon \right) \grad\ub_h(t)}{\grad \vb_h}_{\Elem} = 0
  \end{multline}
  for all $\vb_h\in\VVh^{(n)}_{\Domain}$ and for all $t\in(0,T]$ subject to the initial condition $\ub_h(0) = {\Ext{\ub}}^0$, with $\beta>0$ a penalization factor and $\sigma$ the volume stabilization prefactor.
\end{problem}

Similar to \Cref{sec:tracefem} the discrete covariant derivative $\Grad_{\Elem}$ is described along a componentwise description and extended to the embedding space by using the extended geometric quantities $\nb^{\epsilon}$, $\bm{P}^{\epsilon}$, and $\bm{H}^{\epsilon}$, cf. \cref{eq:gradQu}.

\begin{remark}\label{rem:numericsDI}
  In the considered benchmark problem, see \Cref{sec:NumericalExperiments}, we use an embedding domain $\Domain = [-2,2]^3$ which is discretized by a hierarchical tetrahedral mesh. To be computationally efficient and to ensure a reasonable resolution of $\delta_\epsilon$, an adaptive refinement with about 7--11 grid points across the interface, $\phi_\epsilon \in [0.05,0.95]$, should be used, while a very coarse grid in the remaining part of $\Domain$ is sufficient. For $\Tri{\Domain}$ we define the grid size $h$ by the shortest edge length of the smallest elements, typically located at the interface. To approximate the benchmark surface, we refine the mesh according to the interface thickness of $\epsilon = 0.125$, resulting in a grid size of $h=0.0156$. On this grid we use the \texttt{meshconv} tool \cite{Meshconv2020} to obtain the approximate distance function $\rho_{\epsilon}$. To obtain a numerical approximation with sufficient quality of normals and curvatures as derivatives of $\rho$ requires a proper resolution of the considered surface, see \cite{NV2023diffuse} for a detailed study of these parameters in the vector-valued case. One requirement is that $\epsilon < \delta$, where $\delta$ is the smallest curvature radius of the considered surface. In the benchmark, this implies a very small $\epsilon$, which leads to an unfeasible numerical effort. Therefore, we use the analytic descriptions of $\nb$ and $\bm{H}$ and evaluate and extend them componentwise on $\Tri{\Domain}$ using the Hopf--Lax algorithm. We also consider $\beta = 1000$.
\end{remark}

\subsection{Discussion of the methods}%
\label{sec:Comparison}
We discuss several issues that are important for the numerical treatment of $n$-tensor surface PDEs, in particular the issues listed in the introduction.

First note that there is the following key difference between ISFEM, SFEM, TraceFEM  and DI. The first three methods are based directly on the (variational) PDEs in Problems 1 and 2, which are \emph{consistent} in the sense that they have the same solution, which also coincides with the solution of the $n$-tensor heat equation in the strong formulation. The DI approach, on the other hand, is based on a $\epsilon$-dependent PDE (in a small volumetric neighborhood of the surface), the solution of which, restricted to the surface, is in general different from that of Problems 1 and 2. The formulation only formally converges to the $n$-tensor heat equation as $\epsilon \to 0$.

\emph{Surface representation}.\space
The representation of the surface $\Surface$ is either explicit, in ISFEM and SFEM, or implicit, in TraceFEM and DI. The explicit approach in ISFEM is based on the existence of a parametrization of the surface by an atlas, while SFEM requires only an approximate surface triangulation.
Geometric information (exact or approximate) from the local parametrization at the quadrature points is required for ISFEM quadrature. In SFEM, quadrature is even simpler because only integrals over flat triangles are computed. The implicit description of the surface in TraceFEM and DI is based on a level set description $\phi$ or a phase field description $\phi_\epsilon$ of $\Surface$. In TraceFEM, a surface approximation $\Surface_h$ consisting of triangles is constructed based on a piecewise linear approximation of $\phi$. This requires techniques for computing intersections of tetrahedra with zero levels of linear functions. Due to the fact that the resulting triangulation is in general not shape-regular (``small cuts'') one needs a stabilization (the normal volume derivative stabilization term). As in SFEM, the quadrature is very simple because only integrals over triangles (and tetrahedra) have to be computed. While in TraceFEM an explicit reconstruction $\Surface_h$ of the implicit surface is determined, in the DI method the surface remains implicit. In the discrete variational problems of DI, only integrals over tetrahedra are involved. Thus, quadrature is straightforward. The information of the surface is (only) included via the signed-distance function $\rho$, which is needed in the phase-field function $\phi_\epsilon$. This distance computation requires an additional preprocessing step.

\emph{Representation of the gradient operator and geometry information}.\space
On surfaces there are different natural representations of differential operators of gradient and divergence. In ISFEM the intrinsic representation of the gradient based on local coordinates is used. One then needs a basis of the tangent spaces (at discrete points on the surface).  In ISFEM the orthogonal basis $\{\BasisOrth_1, \BasisOrth_2\}$ is used. The other three methods SFEM, TraceFEM, and DI use a representation of the surface gradient based on the projected standard gradient in $\R^3$.

We now briefly discuss important differences in geometric information between $n=0$ and $n\geq 1$. In the ISFEM, the metric tensor (at discrete points on the surface) is needed for $n\geq 0$, and additionally, for $n \geq 1$, the derivatives of the metric coefficients. For SFEM and TraceFEM, the discrete normal $\nb_h$ is needed for $n\geq0$, while for $n\geq 1$ a more accurate normal approximation (used in $\QQ_h^\sharp$) and an approximation $\bm{H}_h$ of the Weingarten mapping (used in $G_h(\cdot)$) are needed. The DI method requires (approximate) evaluations of the signed distance functions $\rho$ for $n\geq 0$, and additionally (approximate) evaluations of $\nabla \rho$ for $n \geq 1$.

For $n \geq 1$, due to the different representations used, there is the following difference between ISFEM and the other three methods. The methods, SFEM, TraceFEM, and DI, represent the $n$-tensor fields in the embedding space as an element of ${\R^{3^n}}$. For $n \geq 1$, the number of tensor components in the embedding space is larger than in the intrinsic representation used in ISFEM and this discrepancy grows with increasing tensorial rank $n$.

\emph{Tangentiality condition}.\space
Another significant difference between $n=0$ and $n \geq 1$ comes from the tangentiality condition, which is nontrivial only for $n \geq 1$. In ISFEM this condition is automatically satisfied due to the intrinsic representation used.
In SFEM and TraceFEM it is treated by discretizing the augmented variational formulation in \Cref{prob:Problem2}, which includes the consistent penalty term with the projection $\QQ$. This introduces an additional term in the variational form. In the discrete setting, an appropriate scaling of this term is essential. In DI, a volumetric variant $\QQ^\epsilon$ of $\QQ$ is introduced to approximately satisfy the tangentiality condition. Note that on the continuous level, in \Cref{prob:Problem2}, the tangentiality condition is exactly satisfied due to the additional penalty term, while this is not the case for the continuous formulation used in the DI method.

Finally, we will briefly comment on the parameters used in the different methods. In all four methods we have a mesh size parameter $h$, which in ISFEM and SFEM refers to an (approximate) surface triangulation, while in TraceFEM and DI this $h$ corresponds to the mesh size of a tetrahedral triangulation of a volumetric domain containing the surface. In all four methods we have a time step discretization parameter $\Delta t$. In all four methods the polynomial degree $k$ of the finite elements can be chosen. In the presentation above we have restricted ourselves to $k=1$. In ISFEM we have no further parameters. In SFEM and TraceFEM there is a penalty term scaled with $\beta h^{-2}$, so in these methods we have the penalty parameter $\beta$. In TraceFEM, we also have a stabilization term scaled by $\beta' h^{-1}$, so in this method we have the stabilization parameter $\beta'$. In DI there is also a penalty term with a corresponding penalty parameter $\beta$. A key parameter in this method is $\epsilon >0$, which quantifies the interface thickness. The DI method also has a regularization term with a parameter $\sigma$.
The specific parameter values that we use are given below in \Cref{sec:NumericalExperiments}.

\section{Numerical experiments}%
\label{sec:NumericalExperiments}
In this section we present the results of a numerical experiment. We consider an $n$-tensor heat equation, $n=0,1,2$, on a relatively simple surface consisting of a large flat part and a localized bump. The height of this bump is varied and the resulting surfaces have small negative and positive Gaussian curvature values in the bump region (for small bump heights) and (very) large negative and positive Gaussian curvature values in the bump region (for larger bump heights). The initial condition is essentially a regularized Dirac delta function with a support disjoint from the bump support. In \Cref{sec:Experiment} we give a precise description of the problem setting. The four methods described in the sections above are applied to this model problem and some numerical results are presented. The numerical results show that curvature can drastically affect the solution behavior. Specific curvature-related phenomena are discussed in the Sections~\ref{sec:Phenomenon1} -- \ref{sec:Phenomenon3}, for $n=0,1,2$, respectively.

\subsection{Formulation of a tensor diffusion model problem}%
\begin{figure}[ht]
  \centering
  \begin{tikzpicture}
    \node[right] at (0.0,0.0) {\includegraphics{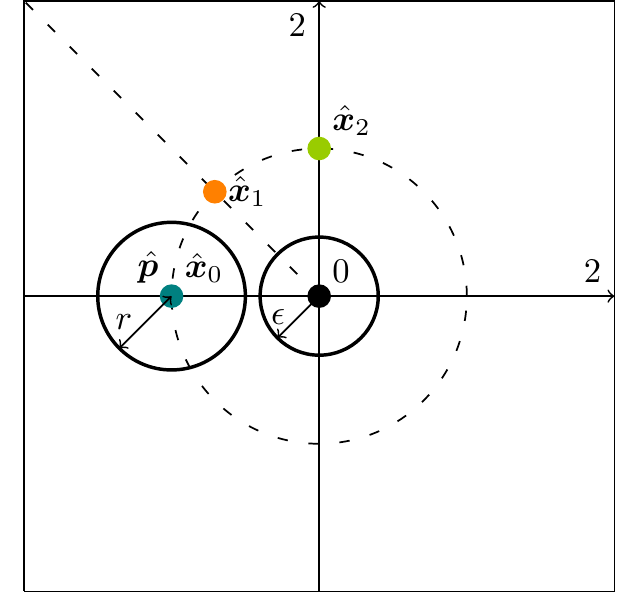}};
    \node[right] at (7.0,-2.0){\includegraphics{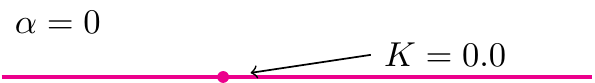}};
    \node[right] at (7.0,0.0) {\includegraphics{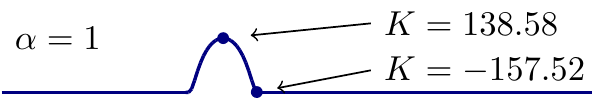}};
    \node[right] at (7.0,2.0) {\includegraphics{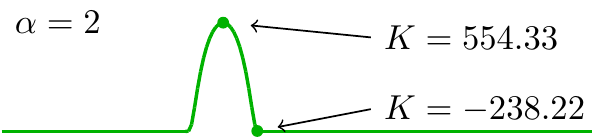}};
  \end{tikzpicture}
  \caption{\label{fig:domain}(Color online) Left: Sketch of the domain with origin colored in black, the outer radius of the bump centered at $\hat{\pb}$ with radius $r$, the initial solution radius $\epsilon$, and the three evaluation points $\xv_0, \xv_1$, $\xv_2$ highlighted in three different colors. The overall domain size of $\RefDomain$ in the numerical computations is chosen to be $[-2,2]^2$. Right: Plot of the bump surfaces along the $\xv_0$-axis for $\alpha\in\{0.0,1.0,2.0\}$. Highlighted are the highest and lowest Gaussian curvature $K$.}
\end{figure}

\label{sec:Experiment}
Let $\Surface$ be the graph of a function $f$, \[\Surface = \Set{ \xb=(\hat{x}^1,\hat{x}^2,f(\hat{x}^1,\hat{x}^2))^T \mid \xv=(\hat{x}^1,\hat{x}^2)^T\in \RefDomain\subset\R^2}\,.\] We want to study a flat surface with an isolated bump that defines a region of negative and positive Gaussian curvature. The bump is described by $f(\xv) \colonequals \alpha \eta(\Norm{\xv - \hat{\pb}}/r)$, where $\alpha \geq 0$ is a scaling factor, $\hat{\pb}\in\RefDomain$ is the center of the bump, and $r > 0$ is its radius. The function $\eta:\R\to\R$ represents a cut-off compressed Gaussian, i.e.,
\[
  \eta(d) = \eta(d; \delta) \colonequals \left\{\begin{array}{ll}
    \exp{\left(-\frac{1}{1 - d^2}\right)} & \text{if }d < 1-\delta \\
    0                                     & \text{otherwise}\,,
    \end{array}\right.
\]
with threshold value $\delta = 0.025$. See \Cref{fig:domain} for a visualization of $\Surface$.

Let $\pb\in\Surface$ be a center point and $\ub_p\in\TensorBundle^n_{\pb}\Surface$ be a (tangential) tensor in $\pb$, then we set as initial condition
\[
  \ub^0(\xb) = \delta_\varepsilon(d_\Surface(\xb, \pb))\,\ub_p\,,\text{ for }\xb\in\Surface\,,
\]
where $\delta_\varepsilon(\cdot)$ is the Dirac delta function. For simplicity, we choose a point $\pb$ in a flat region away from the bump, so that $d_\Surface(\xb, \pb)=\Norm{\xb-\pb}$. The Dirac-delta function is approximated by a single bump of radius $\varepsilon$ around the origin, scaled by $\varepsilon$, so that $\delta_\varepsilon(\Norm{\xb}) = \varepsilon^{-2} \eta(\Norm{\xb}/\varepsilon)$.

We consider $\hat{\pb} = (-0.5, 0.0)^T$ and $r = 0.25$ with varying $\alpha\in[0.0, 2.0]$. For the initial condition we set $\varepsilon=0.2$ and
\[
  \ub_p^{(0)} = 1, \quad \ub_p^{(1)} = (-1,0,0)^T, \quad \ub_p^{(2)} = \ub_p^{(1)}\otimes\ub_p^{(1)}
\]
for the scalar, vector and tensor problem, respectively.

The heat equation is solved in the time interval $t\in[0,1]$ and on the surface $\Surface$ with $\RefDomain = [-2,2]^2$. To illustrate  the solution behavior we define three evaluation points in the parameter domain: $\xv_0=\hat{\pb}$, $\xv_1=0.25\,(-\sqrt{2},\sqrt{2})^T$ and $\xv_2=(0.0,0.5)^T$, all on the circle with radius $0.5$ around the origin in $\RefDomain$, see \Cref{fig:domain}. For the evaluation of the (discrete) solution $\ub_h$, these points have to be lifted to the  discrete surface $\Surface_h$.

The discretization parameters are summarized in \Cref{tab:parameters}.

\begin{table}[ht!]
    \centering
    {\renewcommand{\arraystretch}{1.2}
    \begin{tabular}{c|c|c|c|c|c|c|c}
           & $h$      & $\dt$     & $k$ & $\beta$ & $\beta'$ & $\epsilon$ & $\sigma$ \\
    \hline
   Reference (SFEM)   & $0.0027$  & $10^{-4}$ & 2 & $10$     & ---       & ---        & --- \\
   \hline
  ISFEM    & $0.011$  & $10^{-3}$ & 1   & ---     & ---      & ---        & --- \\
   SFEM    & $0.011$  & $10^{-3}$ & 1   & $10$    & ---      & ---        & --- \\
  TraceFEM & $0.0156$ & $10^{-3}$ & 1   & $0.01$  & $1$      & ---        & --- \\
  DI       & $0.0156$ & $10^{-3}$ & 1   & $10^3$  & ---      & $0.125$    & $10^{-8}$ \\
    \end{tabular}}
    \caption{\label{tab:parameters}Numerical parameters used in the different methods. Note that for the TraceFEM and DI method, the grid size corresponds to the 3d element grid size. The polynomial order $k$ represents the Lagrange polynomial order of the discrete function spaces.}
\end{table}

With the parameters listed in the table, all four methods find approximate solutions within reasonable time on standard hardware.

Previous comparisons of the different methods have shown advantageous properties of SFEM with respect to accuracy and computational effort, see \cite{BJPRV2022Finite}. In order to provide numerical reference data, we use SFEM with a higher spacial and temporal resolution and a higher polynomial order of the solution space. As the finest spacial resolution on the bump, we set the grid size $h\approx0.0027$, the timestep size $\dt=10^{-4}$, and the polynomial degree $k=2$. To reduce the numerical influence of the surface approximation, we have chosen $\ChartMap_h\equiv\ChartMap$ for the SFEM reference computations.

The SFEM and ISFEM methods are implemented with the DUNE/AMDiS framework \cite{Dune,DuneBook,PS2022DuneCurvedGrid}, the TraceFEM method with Netgen/NGSolve and ngsxfem \cite{Schoeberl1997,Schoeberl2014,ngsxfem}, and the DI method within the AMDiS framework \cite{VV2006AMDiS,Witkowskietal_ACM_2015}. The code for the numerical experiments is provided in \cite{Code}.

\subsection{Results for the scalar case}%
\begin{figure}[ht]
  \centering
  \includegraphics{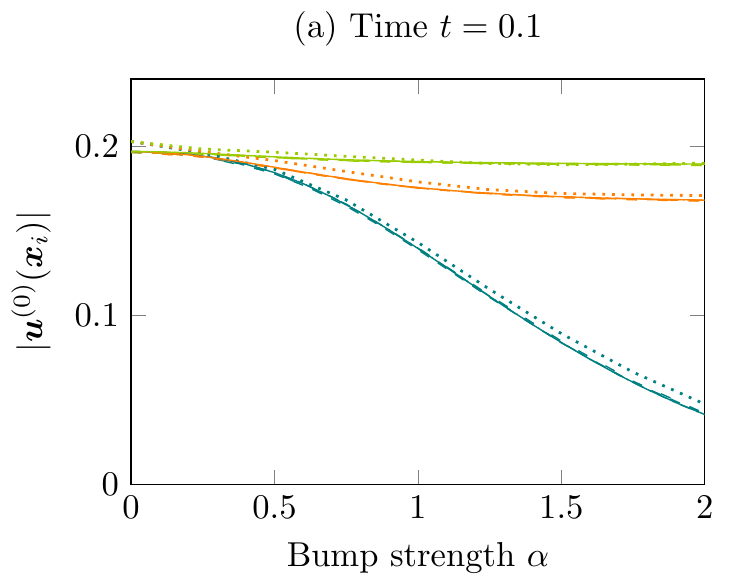}%
  \includegraphics{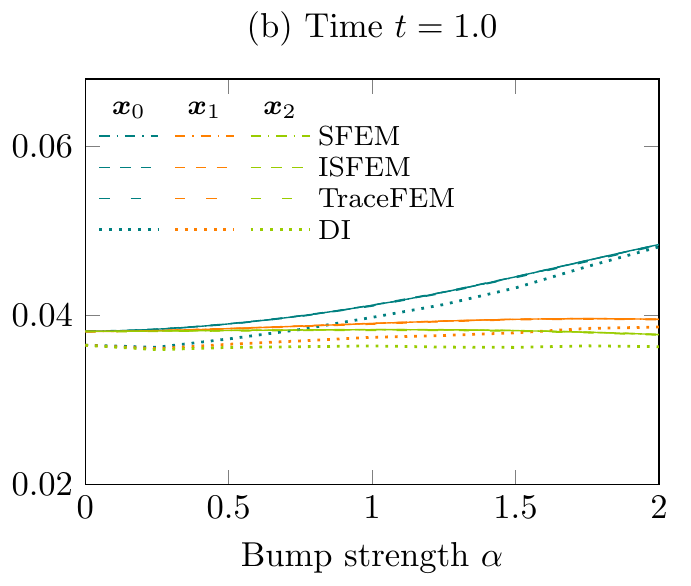} \\
  \includegraphics{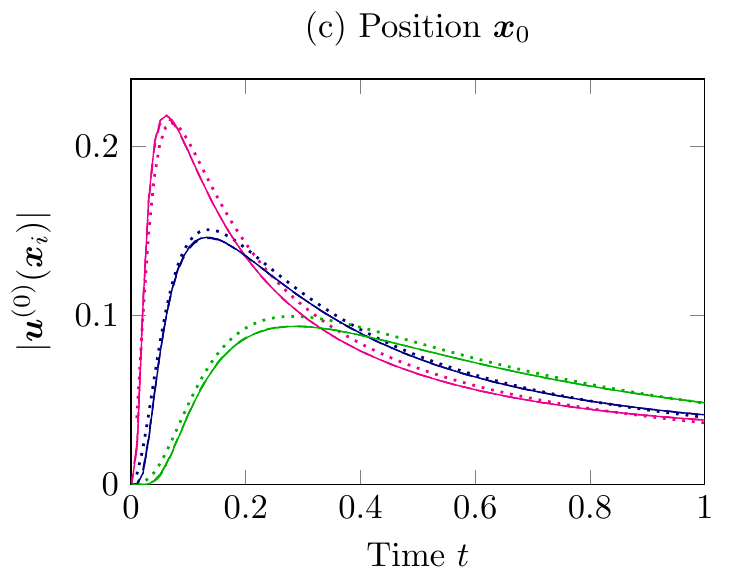}%
  \hspace*{0.2cm}\includegraphics{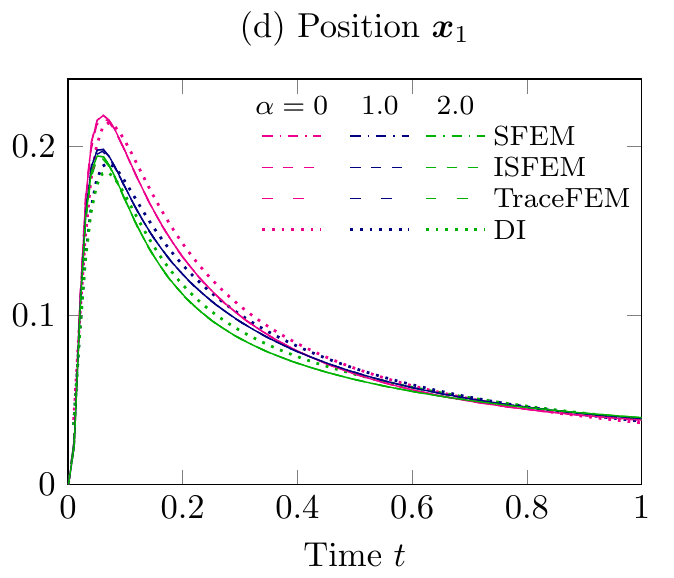}
  \caption{\label{fig:scalar}(Color online) Plot of scalar values $\Abs{\ub^{(0)}(\xb_i)}$ over $\alpha$ (top) and over time $t$ (bottom). Solid lines correspond to a reference solution. Colors correspond to $\xb_0, \xb_1$, and $\xb_2$ and $\alpha\in\{0.0,1.0,2.0\}$.}
\end{figure}

\label{sec:Phenomenon1}
Starting from the initial delta peak at the origin, the scalar heat $\ub^{(0)}$ diffuses over the surface. In flat regions  this diffusion is symmetric. For bump strength $\alpha=0$ this corresponds to the entire domain and thus the maximum heat remains at the initial position. Classical properties can be observed as already described in \Cref{sec:ProblemDefinition}. Not surprisingly, all four methods can represent the flat case equally well, see \Cref{fig:scalar}.

For $\alpha > 0$, the variation of the curvature introduces non-symmetric and anisotropic diffusion into the system. \Cref{fig:scalar} (top) shows that at early times, $t=0.1$, when comparing the solution at the three points, the maximum heat value is at $\xb_2$, while at $t=1.0$ this changes and the maximum value is at $\xb_0$ (on top of the bump). The difference between this maximum value at  $\xb_0$ and the values at the other two points increases for larger $\alpha$ values. Plotted over time in \Cref{fig:scalar} (bottom), there is a transition time point where the maximum changes. This clearly shows that the diffusion depends not only on the geodesic distances on the surface, but also on the surface curvature. The differences at $\xb_1$ and $\xb_2$, points located  in the flat region and having the same distance to the origin, also shows that nearby curved regions influence the solution in the flat part.  On curved surfaces, the simple formula involving the (geodesic) distance holds only for sufficiently short times. An explanation of the phenomenon observed in this experiment can be given by \cref{eq:heat-kernel}. At the bottom of the bump we have \emph{negative} curvature, which leads to \emph{fast} diffusion around the bump, while in a small region containing the bump center, the \emph{positive} curvature \emph{slows down} diffusion, leading to an accumulation of heat in the bump region. The heat diffuses out of the bump region when the difference between the heat values in this region and the region outside the bump is sufficiently large. These local differences also affect nearby zero curvature regions.

The solution behavior is accurately resolved by all four numerical methods. The three consistent methods yield results that (in the ``eye norm'') are hardly distinguishable from the reference solution, while the inconsistent DI method is less accurate (due to a too large $\epsilon$ value).

\subsection{Results for the vector case}%
\label{sec:Phenomenon2}

\begin{figure}[ht]
  \centering
  \includegraphics{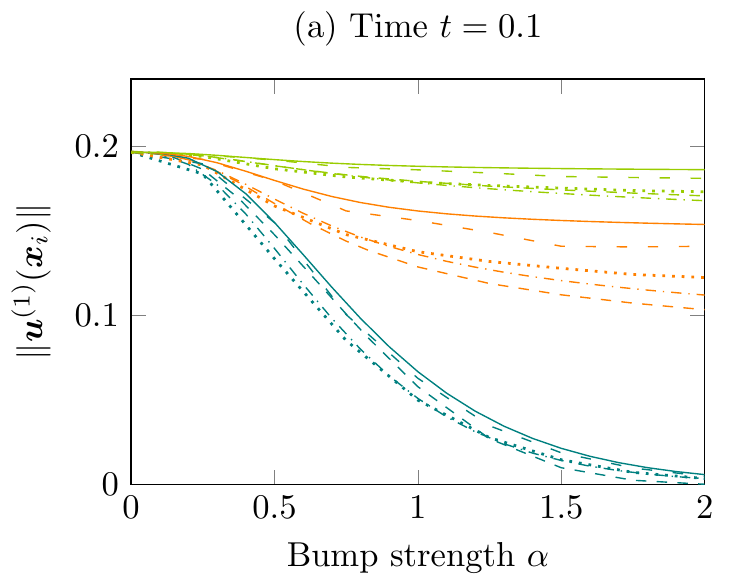}%
  \includegraphics{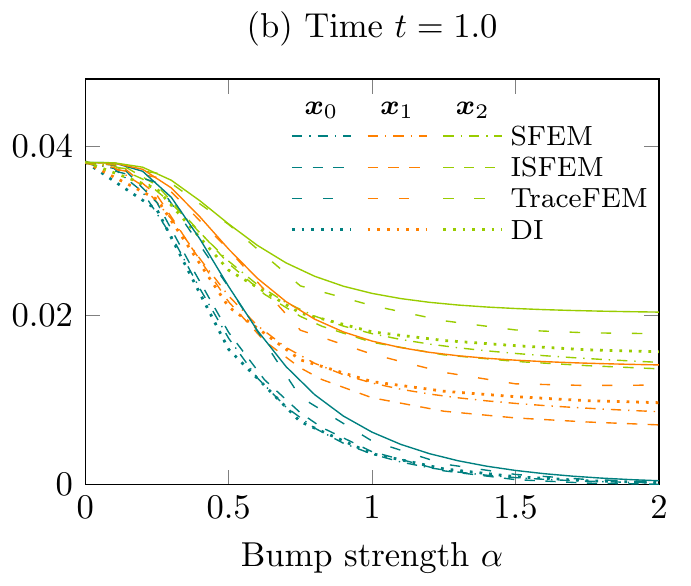}\\
  \includegraphics{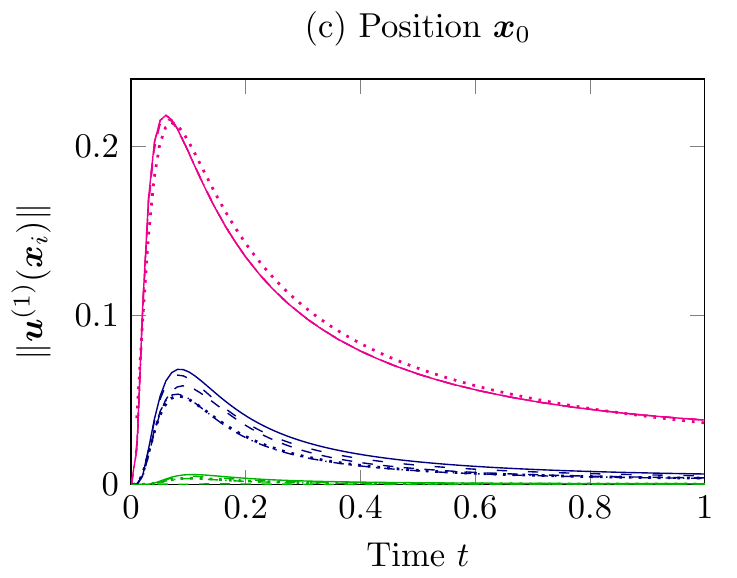}%
  \hspace*{0.2cm}\includegraphics{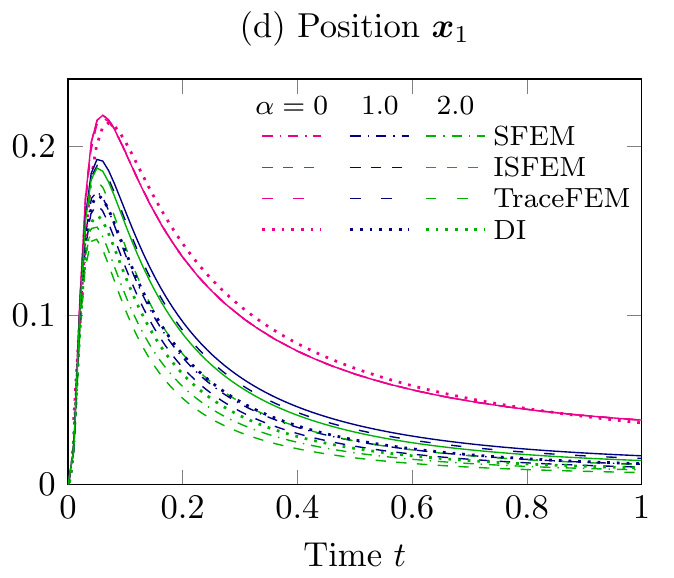}
  \caption{\label{fig:vector-norm}(Color online) Plot of the vector norm $\|\ub^{(1)}(\xb_i)\|$ over $\alpha$ (top) and over time $t$ (bottom). Solid lines correspond to a reference solution. Colors correspond to $\xb_0, \xb_1$, and $\xb_2$ and $\alpha\in\{0.0,1.0,2.0\}$.}
\end{figure}

\begin{figure}[ht!]
  \centering
  \includegraphics{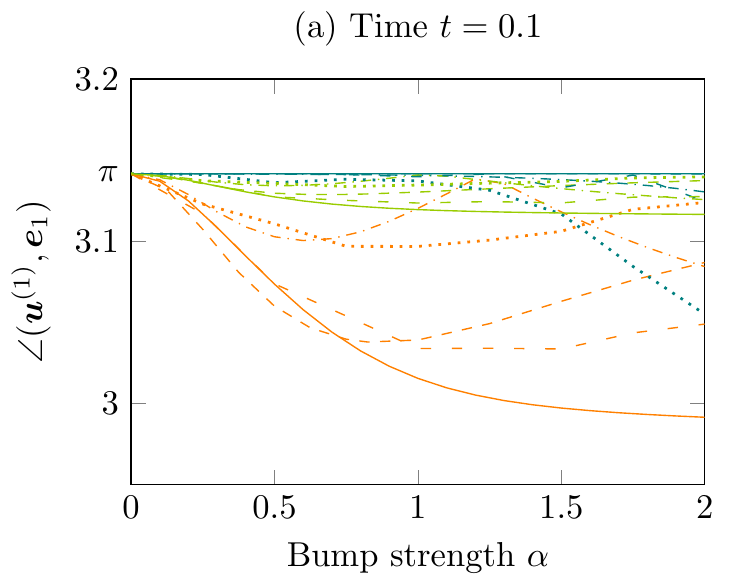}%
  \includegraphics{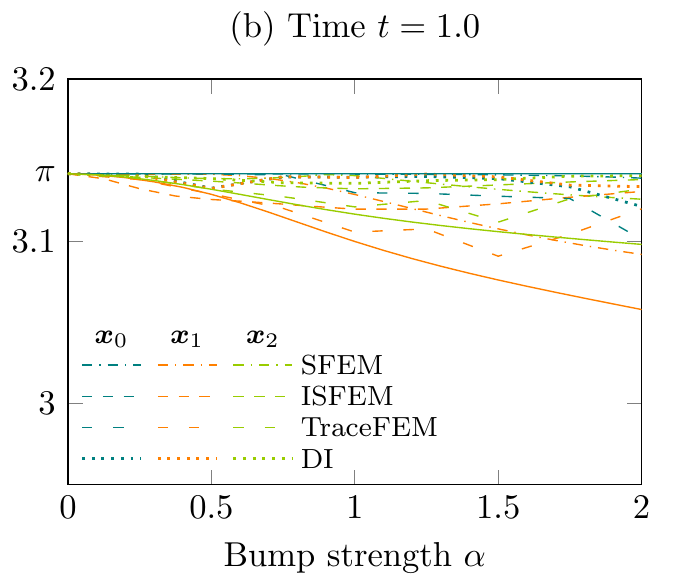}\\
  \includegraphics{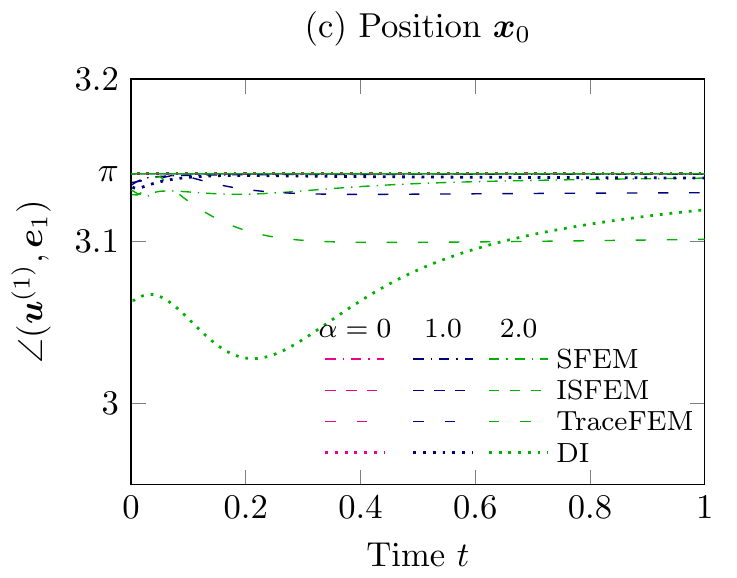}%
  \hspace*{0.2cm}\includegraphics{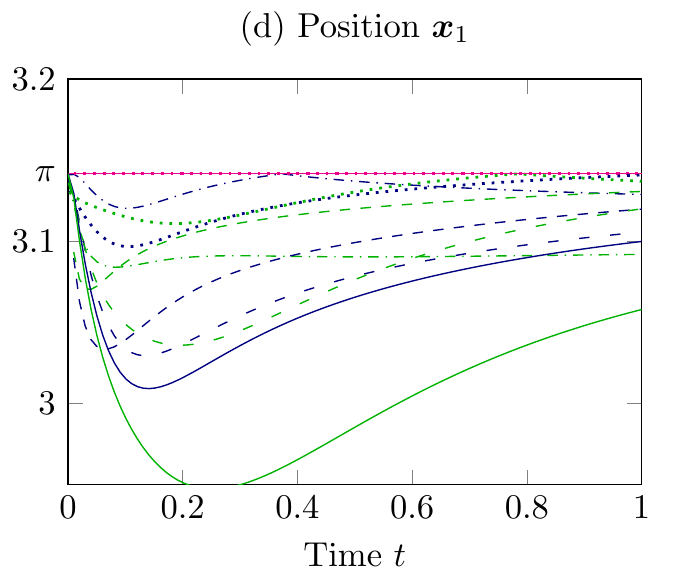}
  \caption{\label{fig:vector-angle}(Color online) Plot of the angle between the vector $\Emb{\ub}^{(1)}(\xb_i)$ and the positive $x$-axis $\bm{e}_1=(1,0,0)^T$ over $\alpha$ (top) and over time $t$ (bottom). Solid lines correspond to a reference solution. Colors correspond to $\xb_0, \xb_1$, and $\xb_2$ and $\alpha\in\{0.0,1.0,2.0\}$.}
\end{figure}

For the vector case not only the norm but also the direction of $\ub^{(1)}$ is of interest. Therefore, we measure the magnitude of the solution $\Norm{\ub^{(1)}}$ and the angle between the vector and the positive $x^1$-axis, i.e., $\angle(\ub^{(1)},\bm{e}_1)\colonequals\arccos\scalarProd{\ub^{(1)}/\Norm{\ub^{(1)}}}{\bm{e}_1}$, at the three reference points, see \Cref{fig:vector-norm} and \Cref{fig:vector-angle}, respectively.

At early times, $t=0.1$, the norm behaves qualitatively similar to the scalar case, but at later times, $t=1.0$, the behavior is very different. While scalar heat diffuses over the whole domain, in the vector case, the norm stays close to zero in the bump center for large $\alpha$, see \Cref{fig:vector-norm}. It does not increase significantly over a very long time. In the scalar case  with $t=1.0$ and $\alpha \in [1,2]$ we see that there is a distinct maximum heat value at the top of the bump, corresponding to $\xb_0$, see \Cref{fig:scalar} (b). In the vector case with $t=1.0$ and $\alpha \in [1,2]$ the opposite happens: the norm values at the top of the bump are much smaller than at the other two points, cf. \Cref{fig:vector-norm} (b).  It seems that there is a strong influence of the additional tangentiality constraint and the interaction with the transport of the direction.

\Cref{fig:vector-angle} shows that the initial \emph{direction} $\ub^0$ is instantaneously extended to the whole domain only in the case of $\alpha=0$. For $\alpha > 0$ this directional extension property holds only near the origin. This is in agreement with the results in \cite{SS2019Vector}, where the limit $t\to 0$ is considered and a vector parallel transport is reconstructed from the vector heat flow solution. For larger times, the curvature of the surface leads to a violation of this property. We see in \Cref{fig:vector-angle} (top) that even in the points in the flat region, i.e., in $\xb_1$ and $\xb_2$ with a flat geodesic to the origin, the ideal angle $\pi$ is missed for $\alpha > 0$. Due to the symmetry of the problem setup, the angle of the solution at the bump, $\xb_0$, is equal to the initial angle. Note that due to the small vector norms on the bump, the evaluation of the angle is poorly conditioned, and thus a small deviation from the $x^1$-axis will result in large deviations in the evaluated angle.

Interpreted as a minimization problem of the Dirichlet energy, $\int_{\Surface}\Norm{\Grad_{\Surface}\ub}^2\,\textrm{d}\xb\to\text{min}$, the vector heat equation minimizes gradients in the magnitude and gradients in the angle. For strongly curved domains, the violation of the angle (caused by the curvature) is compensated by reducing the norm of the vector. This has consequences and leads to increased differences in the norm in the three reference points compared to the scalar case.

All methods show qualitatively the same behavior. However, for all methods the differences to the reference solution are (significantly) larger compared to the scalar case. These differences increase for larger $\alpha$ values. This loss of accuracy compared to the scalar case is caused by the significantly higher numerical complexity for $n\geq 1$, see discussion in \Cref{sec:Comparison}. Depending on the method, the vector case requires the evaluation of derivatives of the projection (SFEM, TraceFEM, DI) or derivatives of the metric coefficients (ISFEM) and is thus become more sensitive to the approximation of the geometry. This is also seen for the vector angle with large variations. These large variations are also due to the low accuracy of the evaluation at the point $\xb_0$ and the sensitivity of the evaluation to small perturbations.

\subsection{Results for the tensor case}%
\begin{figure}[ht]
  \centering
  \includegraphics{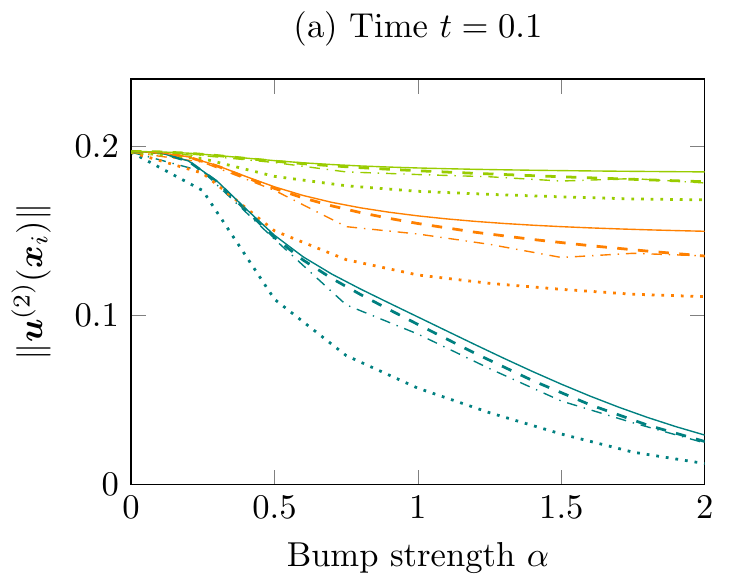}%
  \includegraphics{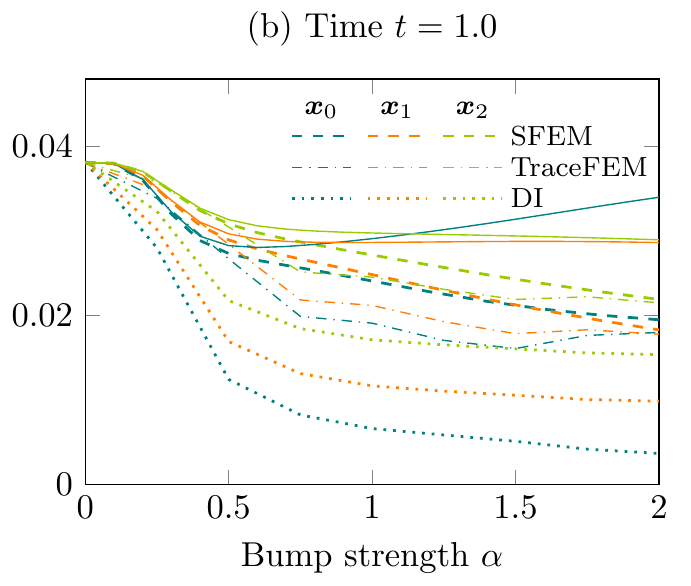}\\
  \includegraphics{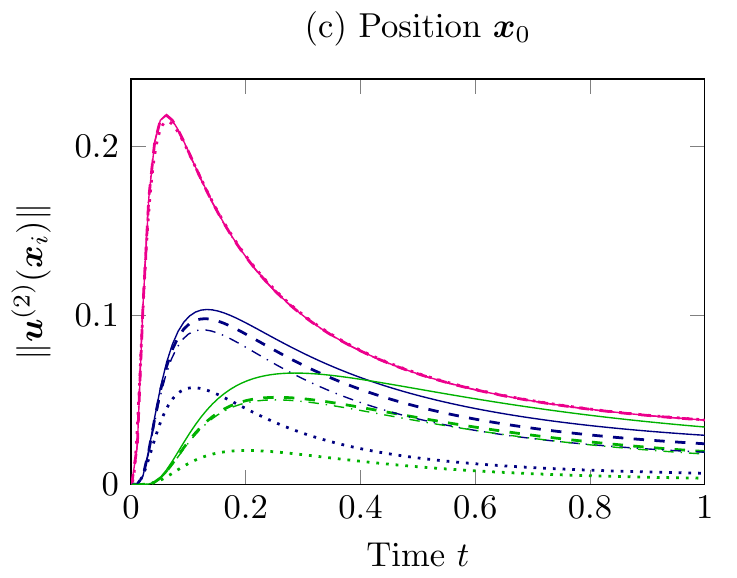}%
  \hspace*{0,2cm}\includegraphics{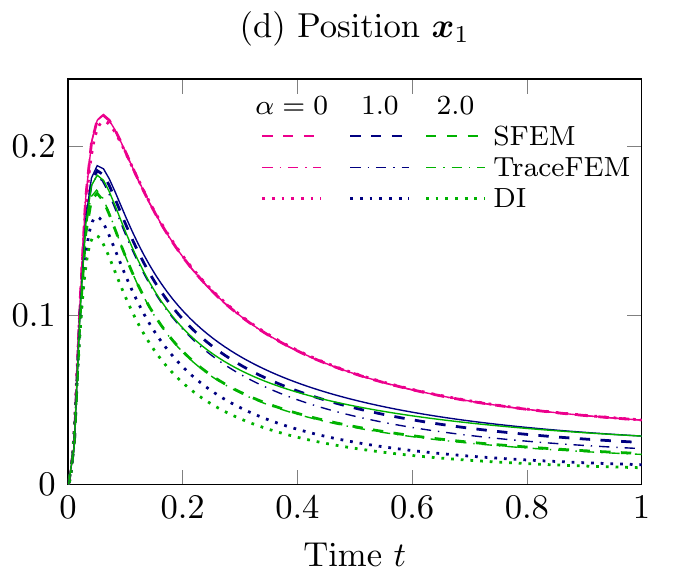}
  \caption{\label{fig:tensor-norm}(Color online) Plot of the tensor norm $\Norm{\Emb{\ub}^{(2)}(\xb_i)}$ over $\alpha$ (top) and over time $t$ (bottom). Solid lines correspond to a reference solution. Colors correspond to $\xb_0, \xb_1$, and $\xb_2$ and $\alpha\in\{0.0,1.0,2.0\}$.}
\end{figure}

\begin{figure}[ht]
  \centering
  \includegraphics{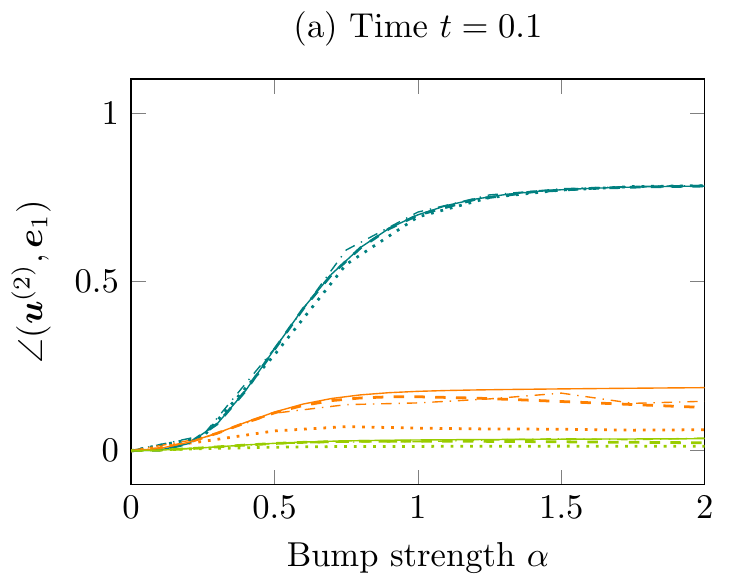}%
  \includegraphics{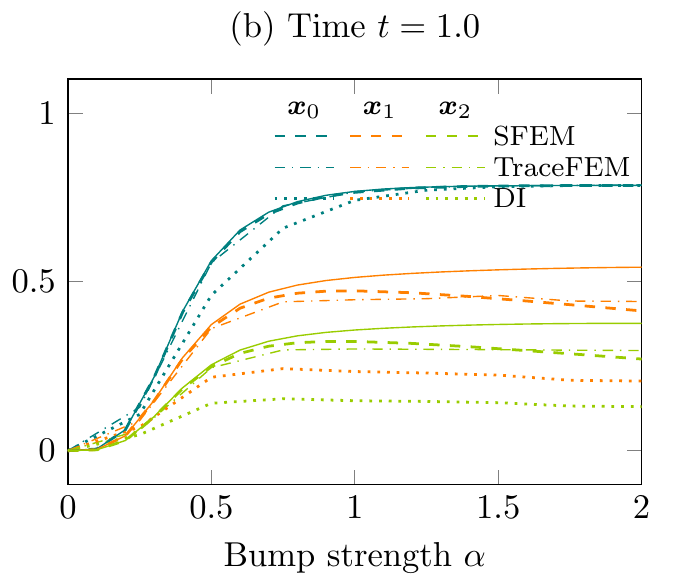}\\
  \includegraphics{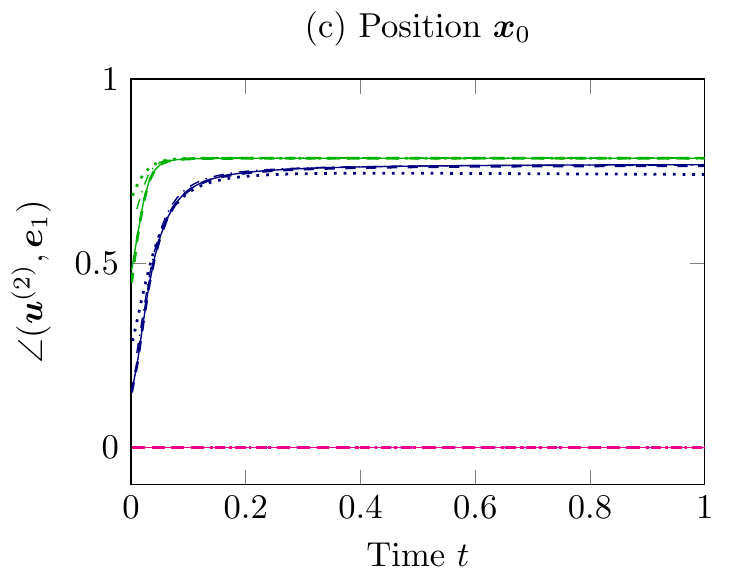}%
  \hspace*{0.2cm}\includegraphics{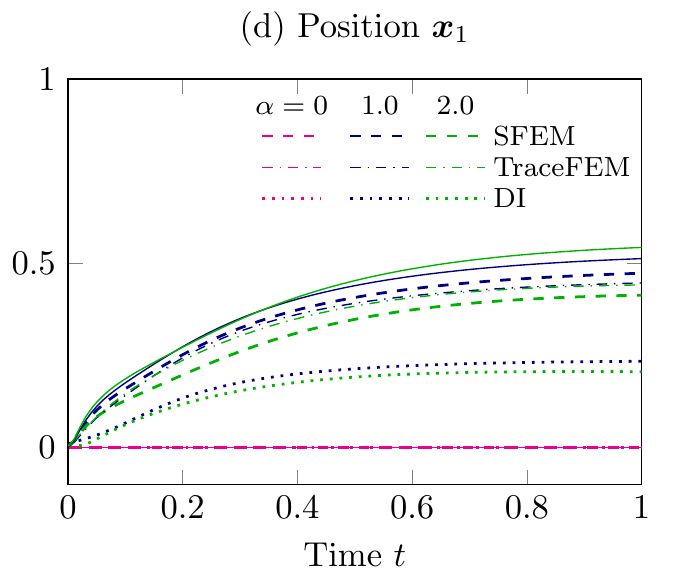}
  \caption{\label{fig:tensor-angle}(Color online) Plot of the angle between the tensor $\ub^{(2)}(\xb_i)$ and the positive $x$-axis $\bm{e}_1=(1,0,0)^T$ over $\alpha$ (top) and over time $t$ (bottom). Solid lines correspond to a reference solution. Colors correspond to $\xb_0, \xb_1$, and $\xb_2$ and $\alpha\in\{0.0,1.0,2.0\}$.}
\end{figure}

\label{sec:Phenomenon3}
For the tensor case we again consider the norm $\Norm{\ub^{(2)}}$ and the angle with the positive $x^1$-axis. The tensor angle is defined as follows: $\angle(\ub^{(2)},\bm{e}_1)\colonequals\arccos\scalarProd{\ub^{(2)}/\Norm{\ub^{(2)}}}{\bm{e}_1\otimes\bm{e}_1}$. Again, we measure these quantities at the three reference points. Due to the increased complexity we here only show results for SFEM, TraceFEM, and DI. The corresponding reference solution is computed on a fine grid, the same as for the scalar and vector case, but with timestep width $\dt=10^{-3}$, see \Cref{fig:tensor-norm} and \Cref{fig:tensor-angle}, respectively. The results are qualitatively similar to the vector case and can be explained by the same reasoning.

The qualitative behavior can be largely resolved by all three methods. However, the differences between the methods continue to increase.

\subsection{Summary}%
\label{sec:Summary}
While some analytical results exist for the diffusion of tangential tensor fields, see \Cref{sec:ProblemDefinition}, quantitative results allowing to test numerical algorithms on simple benchmark problems were missing. We have provided such a setup here. We considered four different numerical methods, ISFEM, SFEM, TraceFEM, and DI, all based on finite element discretizations. They are briefly described and compared. The methods differ with respect to the surface representation, the representation of the gradient operator and geometric information, and the tangentiality condition. The methods are applied to a benchmark problem with a relatively simple surface geometry. We observe that for not too small curvature values the solution behavior is strongly influenced by the geometry. Furthermore, the results show a stronger coupling with geometric properties and an increased sensitivity to the resolution of these properties as the tensor degree increases. Due to this, there is a significant increase in numerical complexity when going from tensor degree $n=0$ to $n \geq 1$.

There are many applications in materials science and biology that exploit the influence of curvature in thin structures. The modeling of such effects often requires tangential vector or tensor fields. We suggest to first test numerical methods for such applications on the provided setup to ensure a proper resolution of the geometric influence.

\paragraph{Acknowledgment}
  The authors wish to thank the German Research Foundation (DFG) for financial support within the Research Unit ``Vector- and Tensor-Valued Surface PDEs'' (FOR 3013) with project no. RE 1461/11-1 and VO 899/22-1. We further acknowledge computing resources provided by ZIH at TU Dresden and within project PFAMDIS at FZ J{\"u}lich.

\bibliographystyle{abbrv_title}
\bibliography{references}
\end{document}